\documentclass{article}

\usepackage{arxiv}

\usepackage[utf8]{inputenc} % allow utf-8 input
\usepackage[T1]{fontenc}    % use 8-bit T1 fonts
\usepackage{hyperref}       % hyperlinks
\usepackage{url}            % simple URL typesetting
\usepackage{booktabs}       % professional-quality tables
\usepackage{amsfonts}       % blackboard math symbols
\usepackage{nicefrac}       % compact symbols for 1/2, etc.
\usepackage{microtype}      % microtypography
\usepackage{lipsum}
\usepackage{graphicx}
\usepackage{amsmath}
\usepackage{caption}
\usepackage{subcaption}
\usepackage{amssymb}
\usepackage{algorithm}
\graphicspath{ {./images/} }

\newtheorem{problem}{Problem}
\newtheorem{lemma}{Lemma}

\newtheorem{remark}{Remark}
\newtheorem{theorem}{Theorem}
\newtheorem{definition}{Definition}

\title{Time-Optimal  Guidance  to Intercept  Moving Targets by Dubins Vehicles}

\author{
 Yuan Zheng \\
  School of Aeronautics and Astronautics\\
  Zhejiang University\\
  Hangzhou 310027, Zhejiang, China\\
   \And
  Xueming Shao \\
   School of Aeronautics and Astronautics\\
  Zhejiang University\\
  Hangzhou 310027, Zhejiang, China\\
  \And
  Zheng Chen* \\
   School of Aeronautics and Astronautics\\
  Zhejiang University\\
  Hangzhou 310027, Zhejiang, China\\
  Corresponding author, Tel. +86-571-87953045\\
  \texttt{z-chen@zju.edu.cn} \\
  \And
  Wenjie Zhao \\
   School of Aeronautics and Astronautics\\
  Zhejiang University\\
  Hangzhou 310027, Zhejiang, China\\
}

\begin{document}
\maketitle
\begin{abstract}
This paper is concerned with a Minimum-Time Intercept Problem (MTIP), for which a Dubins vehicle is guided from a position with a prescribed initial orientation angle to intercept a moving  target in minimum time. Some geometric properties for the solution of the MTIP are presented, showing that the solution path must lie in a sufficient family of 4 candidates. In addition, necessary and sufficient conditions for optimality of each candidate are established. When the target's velocity is constant, by employing the geometric properties,  those 4 candidates are transformed to a class of sufficiently smooth real-valued functions. In order to compute all the 4 candidates, an efficient and robust algorithm to find all the zeros of sufficiently smooth real-valued functions is developed. Since the MTIP with a constant target's velocity is equivalent to the path planning problem of Dubins vehicle in a constant drift field, developing such an algorithm also enables efficiently finding the shortest Dubins path in a constant drift field. Finally, some numerical examples are presented, demonstrating and verifying the developments of the paper.
\end{abstract}

% keywords can be removed

\section{Introduction}

Autonomously guiding  a pursuer to intercept a target in minimum time  is a fundamental problem in the field of guidance and path planning \cite{Isaacs:1965,Merz:1971}.  In this paper, we study a Minimum-Time Intercept Problem (MTIP), for which the target's moving strategy is given and the pursuer is  considered to be a typical nonholonomic vehicle which   moves only forward at a constant speed with  a minimum  turning radius.  Such a nonholonomic vehicle has been commonly dubbed Dubins vehicle in the literature.  As the Dubins vehicle provides an ideal kinematic model for a large class of vehicles, such as fixed-wing unmanned aerial vehicles, autonomous underwater vehicles, unmanned ground vehicles, etc., the shortest paths of Dubins vehicle from a fixed initial configuration (a location and a heading orientation angle) to intercept a target have been widely studied in many fields \cite{Matveev:2012,Matveev:11}.
It should be noted  that for a Dubins vehicle  the shortest path is equivalent to the minimum-time path as the speed is constant.

% the Minimum-Time Interception of a Dubins Vehicle (MTIDV for abbreviation hereafter)   with a target   has been widely studied in many fields.   Therefore, the solution of MTIDV is called either the shortest Dubins path or the minimum-time Dubins path in the literature.
Assuming that the target is stationary and considering that the final impact angle is fixed, the MTIP degenerates to the well-known Dubins problem between two configurations. L. E. Dubins used geometric arguments in \cite{Dubins:57} to show that the shortest Dubins paths between two configurations lie in a sufficient family of $6$ candidates. By relaxing the constraint on the final impact angle, the shortest Dubins path from a configuration to a stationary target was studied in \cite{Boissonnat:1994} and this problem is now called Relaxed Dubins Problem (RDP); the solution of RDP  lies in a sufficient family of 4 candidates.  With the advent of geometric optimal control theory, the developments in \cite{Dubins:57,Boissonnat:1994} were all verified in \cite{Sussmann:94} by using Pontryagin's maximum principle \cite{Pontryagin}. Recently, the shortest path of Dubins vehicle with three consecutive points was studied in \cite{Chen:19automatica} by proposing a polynomial method to compute the solution path.
The syntheses in \cite{Dubins:57,Boissonnat:1994,Sussmann:94,Chen:19automatica} allow computing the shortest Dubins paths within a constant time since only a finite number of candidate paths need to be checked.  However, if the target is moving,  the results developed in  \cite{Dubins:57,Boissonnat:1994,Sussmann:94,Chen:19automatica} do not apply any more for the MTIP.

As a matter of fact, synthesizing the solution of the MTIP with a moving target is fundamentally important in pursuit-evasion engagements \cite{Isaacs:1965}. For this reason, some variants of the MTIP have been studied in the literature. Meyer, Isaiah, and Shima are probably the first ones studying the solution of the MTIP \cite{Mayer:2015}.  Those authors established some sufficient conditions to ensure that the solution of the MTIP is the same as that of the RDP. Whereas, it is not clear what the solution of the MTIP is if the sufficient conditions are not met.
Fixing the final impact angle, the solution of the  MTIP was studied under an assumption that the distance between the initial position and the moving target kept at least 4 times longer than the minimum turning radius \cite{Gopalan:2016}. This assumption restricted the solution into a family of 4 simple candidate paths, allowing formulating some nonlinear equations  so that the optimal path was related to the roots of the nonlinear equations. Those authors proposed using a Newton-type method or a bisection method to find the roots of the nonlinear equations. However, the two numerical methods may not find the desired roots, as shown by the numerical examples in Section \ref{SE:Examples}. More recently, considering that the target moves along a circle, the MTIP was studied in \cite{Manyam:2019,Park:2020}.  In both \cite{Manyam:2019} and \cite{Park:2020}, a strict condition that the initial point of Dubins vehicle is at least 4 times minimum turning radius apart from the target circle was assumed to hold. This strict assumption enables using geometric arguments to synthesize the solution path.

In all the papers cited in the previous paragraph, the target's velocity can be changing. If the target's velocity is constant, it can be proven by a simple coordinate transformation that the MTIP is equivalent to planning the shortest Dubins paths from a configuration to a point in a constant drift field \cite{Techy:2009,Bakolas:2013,McGee:2007}. For this reason, studying the MTIP with a constant target's velocity is quite important in practical scenarios since the motions of fixed-wing unmanned aerial vehicles and autonomous underwater vehicles are usually affected by wind and ocean current, respectively.

Because of the equivalence between the MTIP with a constant target's velocity and the path planning problem in a constant drift field, the MTIP with a constant target's velocity was usually studied from the perspective of synthesizing the shortest Dubins paths in a constant drift field.
McGee, Spry, and Hedrick in \cite{McGee:2005} synthesized the shortest Dubins path in wind and an iterative method was proposed to find the shortest Dubins path. Later on, those authors studied the same problem, showing that the solution lies in a sufficiently family of 8 types \cite[Theorem 1]{McGee:2007}. Analogously, a single implicit equation for the minimum intercept time was formulated and an iterative algorithm was proposed to find the optimal solution in \cite{Looker:2008}, where it was assumed that the distance between the pursuer and the target was always greater than 4 times the minimum turning radius. In all \cite{Looker:2008}, \cite{McGee:2005}, and \cite{McGee:2007}, it was proposed to compute the shortest Dubins paths in constant drift field by iteratively finding a specific zero of nonlinear equations. However, as stated in the previous paragraph, it is challenging to find a specific zero of a nonlinear equation. The reasons include that (1) iterative methods may not converge, and (2) a nonlinear equation usually has multiple zeros so that a zero found by an iterative method cannot be guaranteed to be the desired zero, as shown by the numerical examples in Section \ref{SE:Examples}. In addition to \cite{McGee:2005,McGee:2007}, the properties of shortest Dubins paths in a constant drift field were synthesized in \cite{Bakolas:2013} by using standard optimal control tools and means of discontinuous mapping.  More recently,  the controllability of Dubins problem in a constant drift field was studied in \cite{Caillau:2019}.

Although some variants of the MTIP have been studied in the literature, it is not exaggerate to say that this fundamental problem has not been well addressed, as shown by the counter example presented in \cite[Example 1]{Mayer:2015}. From  practical point of view, it is often required to compute the solution path of the MTIP in real time or onboard, especially for scenarios when the control decisions have to be made in situ or, if not exactly, at least efficiently. However, to the authors' best knowledge, an efficient and robust algorithm for computing the solution of the MTIP does not exist in the literature.

Unlike the aforementioned papers, the solution of the MTIP is thoroughly investigated in this paper without any assumption on the distance between the initial point and the target. First, by introducing three functions in terms of the target's trajectory and studying the continuity properties of the three functions, some geometric properties of the solution path are presented. Using these geometric properties, it is proven that the solution path of the MTIP lies in a sufficient family of 4 candidates, and necessary and sufficient conditions for optimality of each candidate are established. In addition, when considering that the target's velocity is constant, these geometric properties are used  to formulate some nonlinear equations so that the solution of the MTIP is determined by a specific zero of the nonlinear equations.

In general, a nonlinear equation may have multiple zeros but only one specific zero is related to the solution of the MTIP. Since the typical Newton-like iteration method and bisection method proposed in \cite{McGee:2005,McGee:2007} cannot be guaranteed to find a specific zero of a nonlinear equation,  a new algorithm is developed in this paper to find all the zeros of a sufficiently smooth real-valued function. Applying this algorithm allows computing the solution of the MTIP within a constant time, which provides a potential for onboard applications. It is worth mentioning that developing such an algorithm also enables computing the shortest Dubins path in a constant drift field, as shown by the last numerical example in Section \ref{SE:Examples}.

 The paper is organized as follows. In Section \ref{SE:Preliminary} the MTIP is formulated as an optimal control problem and its necessary conditions for optimality are presented according to Pontryagin's maximum principle. Section \ref{SE:Syntheses} is attributed to establishing geometric properties for the solution of the MTIP. Then, a robust and efficient algorithm is designed in Section \ref{SE:Algorithm} to find the solution of the MTIP. All the developments are finally demonstrated and verified by numerical examples in Section \ref{SE:Examples}.

\section{Preliminary}\label{SE:Preliminary}

In this section, the MTIP is formulated as an optimal control problem, and its necessary conditions for optimality are established according to Pontryagin's maximum principle.

\subsection{Problem Formulation}

Consider a 2-dimensional engagement involving a pursuer and a moving  target. The pursuer is a Dubins vehicle that moves only forward at a constant speed with a minimum turning radius. Denote the state (or configuration) of the pursuer by $\boldsymbol{z}:=(x,y,\theta)\in \mathbb{R}^2\times \mathbb{S}^1$, which consists of a position vector $(x,y)\in \mathbb{R}^2$ and a heading orientation angle $\theta \in \mathbb{S}^1$. Then, by normalizing the position $(x,y)$ so that the pursuer's speed is one, the kinematics is expressed as
\begin{align}
\qquad \dot{\boldsymbol{z}}(t)=\left[
\begin{matrix}
 \cos \theta(t)\\
 \sin \theta(t)\\
 u(t)/\rho
\end{matrix}
\right],\ \ \ \ \ \ \ \ \ \ \  u\in [-1,1]
\label{Eq:problem1}
\end{align}
where $t\geq 0$ denotes time, the dot denotes the differentiation with respect to time,  $\rho>0$ is the minimum turning radius,  and $u$ is the control input representing the lateral acceleration of the pursuer.  Without loss of generality, we assume that the configuration at initial time $t=0$ is
$$\boldsymbol{z}_0 := (0,0,\pi/2).$$

%Since the target is nonmaneuvering, it follows that the velocity  $\boldsymbol{v}$ is a constant vector.

 Denote by $\boldsymbol{v} = (v_x,v_y)\in \mathbb{R}^2$  the velocity of the target.  Let the position of the target at initial time $t=0$ be $(\hat{x}_{0},\hat{y}_{0})$. Then, the position of the target at any time $t\geq 0$ is given by
\begin{align}
E(t) = (\hat{x}_{0},\hat{y}_{0})+ \int_{0}^t \boldsymbol{v}(\tau) \mathrm{d}\tau.
\label{EQ:E(t)}
\end{align}
Throughout the paper, we assume that the target's moving strategy is given so that the position $E(t)$ for any $t\geq 0$ is available to the pursuer.

The MTIP is an optimal control problem defined as below.
\begin{problem}[MTIP]\label{problem1}
The MTIP consists of finding a minimum time $t_m > 0$ so that  the system in Eq.~(\ref{Eq:problem1}) is steered by a measurable control $u(\cdot)$ over the interval $[0,t_m]$ from a fixed initial configuration $\boldsymbol{z}_0$ at $t=0$ to intercept a moving target at $t_m$, i.e., $(x(t_m),y(t_m)) = E(t_m)$.
\end{problem}
If $\boldsymbol{v} \equiv 0$, the MTIP degenerates to the well-known RDP \cite{Boissonnat:1994}. If the velocity $\boldsymbol{v}$ is constant, by a simple coordinate transformation, it can be proven that the MTIP is equivalent to planning path for a Dubins vehicle from a fixed configuration to a fixed point \cite{Bakolas:2013}.

 In \cite{Cockayne:1967}, it was shown that the solution of the MTIP exists when $ \|{\boldsymbol{v}}\|< 1$, where the notation $\|\cdot\|$ denotes the Euclidean norm. In this paper, we assume that the condition $\|{\boldsymbol{v}}\|<1$ stands  so that the solution of the  MTIP exists.

\subsection{Necessary Conditions}

Let  $p_x,\ p_y$ and $p_{\theta}$ be the costate variables of $x,\ y$, and $\theta$, respectively. Then, the Hamiltonian of the MTIP  is
\begin{align}
H=p_x \cos \theta +p_y \sin \theta +p_{\theta}u/\rho-1
\label{Eq:problem4}
\end{align}
According to Pontryagin's maximum principle \cite{Pontryagin}, we have
\begin{align}
%\begin{cases}
\dot{p}_x(t)&=-\frac{\partial H}{\partial x}=0\label{EQ:px}\\
\dot{p}_y(t)&=-\frac{\partial H}{\partial y}=0\label{EQ:py}\\
\dot{p}_{\theta}(t)&=-\frac{\partial H}{\partial \theta}=p_x(t)\sin \theta(t)-p_y(t) \cos \theta(t)
%\end{cases}
\label{Eq:problem5}
\end{align}
It is apparent from Eq.~\eqref{EQ:px} and Eq.~\eqref{EQ:py} that $p_x$ and $p_y$ are constant. By integrating  Eq.~(\ref{Eq:problem5}), we have
\begin{align}
p_{\theta}(t)=p_xy(t)-p_yx(t)+c_0
\label{Eq:problem6}
\end{align}
where $c_0$ is a constant.  In view of Eq.~(\ref{Eq:problem6}),  if $p_{\theta}\equiv 0$ on a nonzero interval, the path of $(x,y)$ is a straight line segment on this interval. Note that $u\equiv 0$ along a straight line.  Thus,  we have  $u\equiv0$ if $p_{\theta} \equiv 0$.  As a result, the maximum principle indicates that the optimal control $u$ is totally determined by $p_{\theta}$, i.e.,
\begin{align}
u=\begin{cases}
1,&p_{\theta}>0\\
0,&p_{\theta}\equiv0\\
-1,&p_{\theta}<0
\end{cases}
\label{Eq:problem9}
\end{align}
The path $t\mapsto [x(t),y(t)]$ is a circular arc with right (resp. left) turning direction if $u = -1$ (resp. $u=1$). Therefore, the switching conditions in Eq.~\eqref{Eq:problem9} imply that  the solution path of the MTIP is a concatenation of circular arcs and straight line segments.

The necessary conditions from Eq.~\eqref{EQ:px} to Eq.~\eqref{Eq:problem9} will be used in the next section to synthesize the solution path of the MTIP. Before proceeding, some notations and definitions are defined in the following subsection.

\subsection{Notations and definitions}

%Denote by $E:[0,+\infty]\rightarrow \mathbb{R}_+$ the path of the target with respect to time, i.e.,
%\begin{align}
%E(t) = (\hat{x}_0,\hat{y}_0) + \boldsymbol{v} t, \ \ \ t\geq  0.
%\end{align}
%By the definition of $t_m$ in Problem \ref{problem1}, it is clear that $E(t_m)$ is the point at which the target is intercepted by the pursuer with a minimum time.
Denote by
\begin{align}
F:\mathbb{R}^2\rightarrow [0,+\infty),\ \  (x,y)\mapsto F[x,y]
\end{align}
 the minimum time (equivalent to the length of the shortest path since the speed of the pursuer is $1$) for the pursuer to move from the initial configuration $\boldsymbol{z}_0$ to the point $(x,y)\in \mathbb{R}^2$. By definition, it is apparent that the value of $F[x,y]$ denotes the length of the RDP's path from $\boldsymbol{z}_0$ to $(x,y)$. According to \cite{Boissonnat:1994}, the solution of the RDP can be computed in a constant time by checking at most four candidate paths. Thus, given any $t\geq 0$, the value of $F[E(t)]$ is readily available.

 %Notice that  $M(t)$ may not be continuous as will be shown in the subsequent analysis.

\begin{figure}[h]
	\centering
	\includegraphics[width = 7cm]{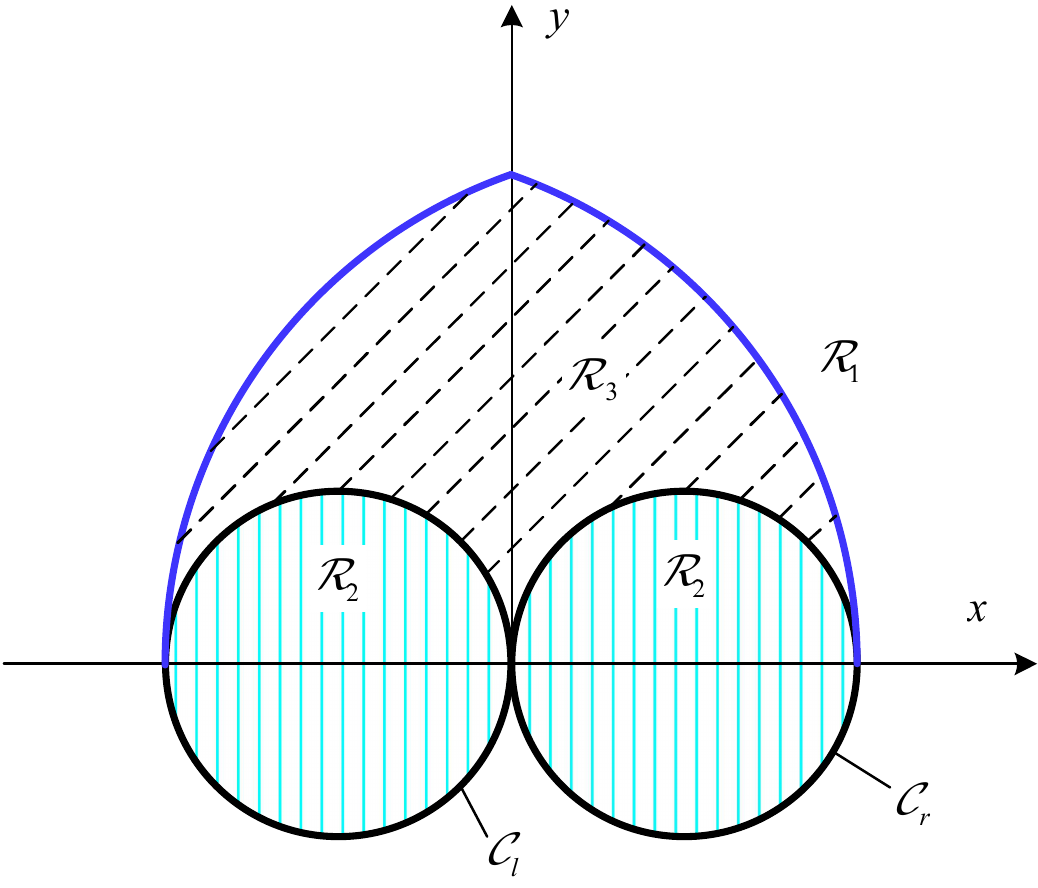}
	\caption{The geometry of the circles $\mathcal{C}_r$ and $\mathcal{C}_l$ and the subregions $\mathcal{R}_1$, $\mathcal{R}_2$, and $\mathcal{R}_3$ in  the 2-dimensional plane.}
	\label{Fig:synthese6}
\end{figure}

Let $\mathcal{C}_r$ and $\mathcal{C}_l$ be circles of radius $\rho$, lying on the right and left side  of  initial configuration $\boldsymbol{z}_0$, respectively, i.e.,
\begin{align}
\mathcal{C}_r = \{(x,y)\in \mathbb{R}^2 | (x-\rho)^2 + y^2 = \rho^2\}\nonumber
\end{align}
and
\begin{align}
 \mathcal{C}_l = \{(x,y)\in \mathbb{R}^2| (x+\rho)^2 + y^2 = \rho^2\}\nonumber
\end{align}
Moreover, set
\begin{align}
\mathcal{D}_r = \{(x,y)\in \mathbb{R}^2 | (x-\rho)^2 + y^2 \leq  \rho^2\} \nonumber
\end{align}
and
\begin{align}
\mathcal{D}_l = \{(x,y)\in \mathbb{R}^2 | (x+\rho)^2 + y^2 \leq \rho^2\}\nonumber
\end{align}
as the regions bounded by $\mathcal{C}_r$ and $\mathcal{C}_l$, respectively.
Denote by
\begin{align}
\boldsymbol{c}_0^{r} :=
\left(
\begin{array}{c}
\rho \\
0
\end{array}
\right) \ \text{and}\
\boldsymbol{c}_0^{l} := \left(
\begin{array}{c}
-\rho \\
0
\end{array}
\right)\nonumber
\end{align}
the centers of $\mathcal{C}_r$ and $\mathcal{C}_l$, respectively.
We define three subregions $\mathcal{R}_1$, $\mathcal{R}_2$, and $\mathcal{R}_3$ in the $2$-dimensional plane as follow:
\begin{align}
\mathcal{R}_2 =&\ \mathcal{D}_r \cup \mathcal{D}_l,\nonumber\\
\mathcal{R}_{3} = &\ \big \{(x,y)\in \mathbb{R}^2 \ |\  y > 0, (x -\rho)^2 + y^2 \leq 9\rho^2, (x+\rho)^2 + y^2 \leq 9 \rho^2\big\} -( \mathcal{R}_2\cap \{y>0\}), \nonumber\\
\mathcal{R}_{1} =&\  \mathbb{R}^2 - \mathcal{R}_2 - \mathcal{R}_3.\nonumber
\end{align}
The geometries for $\mathcal{C}_r$, $\mathcal{C}_l$, $\mathcal{R}_1$, $\mathcal{R}_2$, and $\mathcal{R}_3$ are all illustrated in Fig.~\ref{Fig:synthese6}.

 Denote by ``S'' and ``C'' a straight line segment and a circular arc with radius of $\rho$, respectively. If a circular arc C has a right (resp. left) turning direction, we represent it by R (resp. L).  In addition, we denote by ``S$_d$'' a straight line segment with a length of $d\geq 0$, and denote by ``C$_{\alpha}$'' a circular arc with a radian of $\alpha\geq 0$.  Accordingly, we denote by L$_{\alpha}$ (resp. R$_{\alpha}$) a left-turning (resp. right-turning) circular arc with its radian being $\alpha \geq 0$.
\begin{definition}[Feasible Dubins path]
Given any path in the $x-y$ plane, it is said a feasible Dubins path if the curvature along the path everywhere is not greater than $1/\rho$.
\end{definition}

As the solution of the RDP is used in the following text, we summarize the solution types of the RDP  by the following remark.
\begin{remark}[J.-D. Boissonnat and X.-N. Bui \cite{Boissonnat:1994}]\label{RE:RDP}
The solution path of the RDP belongs to either CS or CC or substrings thereof, where
\begin{itemize}
\item CC = \{RL, LR\}
\item CS = \{RS, LS\}
\end{itemize}
\end{remark}

\section{Characterizing the solution  of the MTIP}\label{SE:Syntheses}

In this section, some geometric properties for the solution path of the MTIP will be established by analyzing the function $F[E(t)]$. For simplicity of presentation, all the proofs for the theorems of this section are postponed to Appendix \ref{Appendix:A}.

By the following lemma, we first recall from \cite{Mayer:2015} the properties for the solution of the MTIP in the case that $F[E(t)]$ is continuous.
\begin{lemma}[Meyer, Isaiah, and Shima \cite{Mayer:2015}]\label{LE:continuous}
If the function $F[E(t)]$ is continuous for $t\in [0,+\infty)$, then the following two statements hold:
\begin{description}
\item (1) The minimum interception time $t_m>0$ is the minimum fixed point of $F[E(t)]$, i.e.,
\begin{align}
t_m = \mathrm{min} \{t > 0\ \big| \  t = F[E(t)]\}.\nonumber
\end{align}
\item (2) The solution path of the MTIP is the same as that of the RDP from the initial condition $\boldsymbol{z}_0$ to the interception point $E(t_m)$.
\end{description}
\end{lemma}

This lemma presents the relationship between $t_m$ and $F[E(t)]$ under the sufficient condition that the function $F[E(\cdot)]$ is continuous. However, this sufficient condition may not be met, as shown by the following lemma.
\begin{lemma}\label{TH:continuity}
The function $F[E(t)]$ is discontinuous at a time $\bar{t}>0$ if and only if
\begin{align}
E(\bar{t}) \in \{(x,y)\in  \mathcal{C}_r \cup \mathcal{C}_l\ \big|\  y>0 \}.\nonumber
\end{align}
\end{lemma}
%Once the function $F[E(t)]$ is discontinuous, the relationship between $t_m$ and $F[E(t)]$ is not known in the literature.
By extending Lemma \ref{LE:continuous}, the following lemma presents the solution property of the  MTIP  without requiring the continuity of $F[E(t)]$.
\begin{lemma}\label{LE:fixed_RDP}
No matter the function $F[E(t)]$ is continuous or not, if the minimum intercept time $t_m>0$ is a fixed point of $F[E(t_m)]$, i.e., $t_m = F[E(t_m)]$, the solution path of the MTIP  is the same as that of the RDP from $\boldsymbol{z}_0$ to $E(t_m)$.
\end{lemma}

%In the sequel, the solution path of MTIP will be investigated by considering the location of $E(t_m)$ with respect to three subregions $\mathcal{R}_1$, $\mathcal{R}_2$, and $\mathcal{R}_3$.

In order to establish necessary and sufficient conditions for $t_m = F[E(t_m)]$, we consider the location of $E(t_m)$ in different subregions of $\mathcal{R}_1$, $\mathcal{R}_2$, and $\mathcal{R}_3$ in the sequel.
\begin{theorem}\label{LE:Occurance_R1}
If the final point $E(t_m)$ of the MTIP lies in the interior of $\mathcal{R}_1\cup \mathcal{R}_2$, i.e., $E(t_m) \in \mathrm{int}(\mathrm{R}_1\cup \mathcal{R}_2)$, then the minimum interception time $t_m$ is a fixed point of $F[E(\cdot)]$, i.e., $t_m = F[E(t_m)]$.
\end{theorem}

Combining Lemma \ref{LE:fixed_RDP} and Theorem \ref{LE:Occurance_R1}, we have that, no matter the function $F[E(t)]$ is continuous or not, only if $E(t_m) \in \mathrm{int}(\mathrm{R}_1\cup \mathrm{R}_2)$, the solution path of the MTIP must be the same as that of the RDP from $\boldsymbol{z}_0$ to $E(t_m)$. Therefore, according to Remark \ref{RE:RDP}, we immediately have the following result.
 \begin{remark}\label{RE:R1_R2}
 If $E(t_m) \in \mathrm{int}(\mathrm{R}_1\cup \mathcal{R}_2)$, the solution path of the MTIP belongs to CC or CS or substrings thereof.
 \end{remark}

Next, we shall establish the geometric properties of the solution of the MTIP for the rest case that $E(t_m) \in \mathrm{int}(\mathcal{R}_3)$. Before proceeding, we present a symmetric property by the following remark.
\begin{remark}[Symmetric property \cite{Boissonnat:1994}]\label{RE:Symmetric}
Given any feasible Dubins path from $\boldsymbol{z}_0$ to $(x,y)\in \mathbb{R}^2$, there exists a feasible Dubins path from $\boldsymbol{z}_0$ to $(-x,y)$ so that the two feasible Dubins paths are symmetric with respect to the $y$ axis.
\end{remark}
%Notice that the solution of MTIP is a feasible Dubins path from $\boldsymbol{z}_0$ to $E(t_m)$.
Thanks to this remark, in the following paragraphs we just consider the solution of the MTIP with the terminal point $E(t_m)$ in the subregion $\mathcal{R}_3\cap \{x> 0\}$.

\begin{figure}[h]
	\centering
	\includegraphics[width = 8cm]{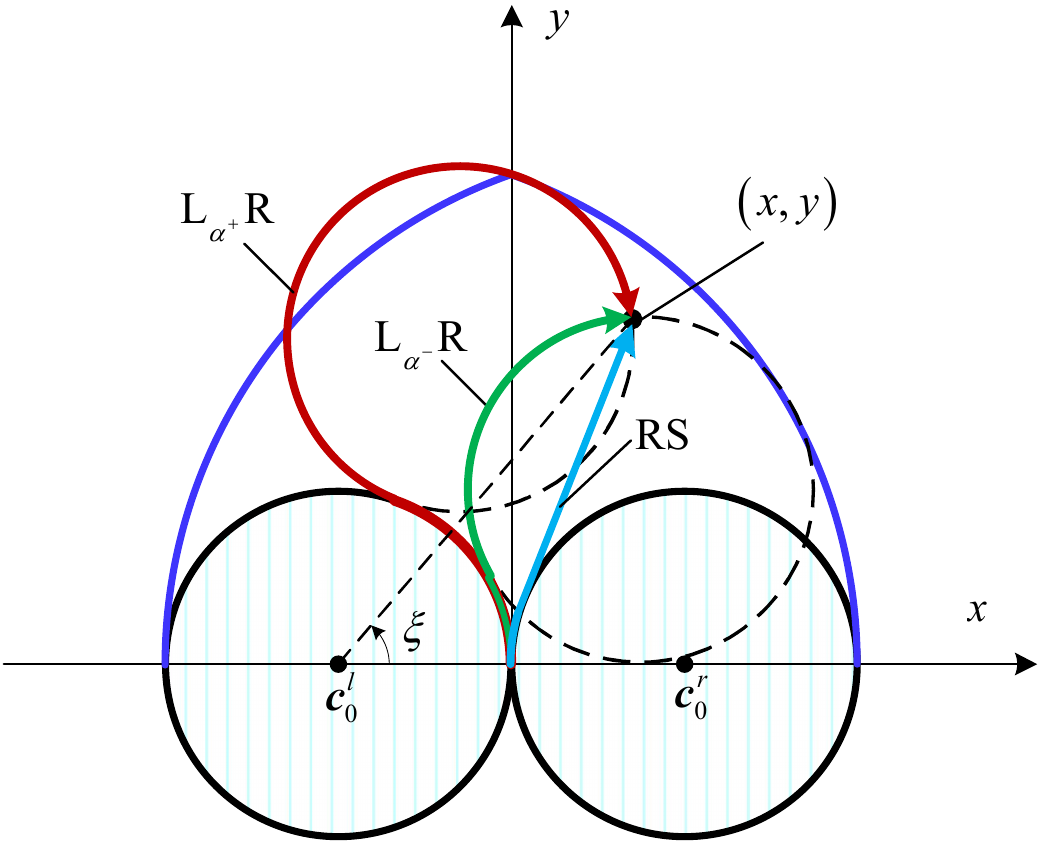}
	\caption{The geometries of the feasible Dubins paths related to $F$, ${L}^-$, and ${L}^+$.}
	\label{Fig:L1L2}
\end{figure}

Notice that for any point $(x,y)\in \mathrm{int}(\mathcal{R}_3\cap \{x>0\})$, there exist two circles of radius $\rho$ that not only pass through $(x,y)$ but also are tangent to $\mathcal{C}_l$, as illustrated in Fig.~\ref{Fig:L1L2}.
Circular arcs on the two circles together with circular arcs on  $\mathcal{C}_l$ form two feasible Dubins paths of type LR from $\boldsymbol{z}_0$ to $(x,y)$, and we denote by L$_{\alpha^+}$R and L$_{\alpha^-}$R the types of the two feasible Dubins paths, as shown in Fig.~\ref{Fig:L1L2}.

 For any $(x,y)\in \mathcal{R}_3\cap \{x> 0\}$, let  $\xi>0$ denote the angle between the positive $x$-axis and the vector from $\boldsymbol{c}_0^l$ to $(x,y)$, as shown in Fig.~\ref{Fig:L1L2}.  It is apparent from Fig.~\ref{Fig:L1L2} that $\alpha^+ \geq \xi$ and $\alpha^- \leq \xi$ \cite{Ding:2019}.
 \begin{definition}
 %For any point $(x,y)\in \mathcal{R}_3\cap \{x> 0\}$, we denote by $L_{\alpha\leq\xi}R$ and $L_{\alpha\geq\xi}R$ the types of the feasible Dubins paths from $\boldsymbol{z}_0$ to $(x,y)$, as illustrated in Fig.~\ref{Fig:L1L2}. And,
 We denote by  $L^+[x,y]$ and $L^-[x,y]$  the lengths of the feasible Dubins paths from $\boldsymbol{z}_0$ to $(x,y)\in \mathcal{R}_3\cap\{x>0\}$ with types $L_{\alpha^+}R$ and $L_{\alpha^-}R$, respectively.
 \end{definition}
 The following lemma shows the continuity properties of $L^+[E(t)]$ and $L^-[E(t)]$.
  %for the time interval such that $E(t) \in \mathcal{R}_3\cap \{x>0\}$.

\begin{lemma}\label{LE:continuity}
Let $\bar{t}>0 $ be a time so that $E(\bar{t}) \in \mathcal{R}_3\cap \{x>0\}$. Then, both $L^-[E(t)]$ and $L^+[E(t)]$ are continuous at $\bar{t}$.
%For any continuous curve $\gamma$ in $\mathcal{R}_3\cap\{x>0\}$, the two functions $L^-(\gamma)$ and $L^+(\gamma)$ are continuous.
\end{lemma}

%Since the proof of this lemma involves some basic mathematic analysis, it  is delayed to Appendix \ref{Appendix:A}.

We  recall from \cite{Ding:2019} the relationship of the three functions $F$, $L^-$, and $L^+$ by the following lemma.
\begin{lemma}[Ding, Xin, and Chen \cite{Ding:2019}]\label{LE:Ding}
Given any point $(x,y)$ in $\mathrm{int}(\mathcal{R}_3\cap \{x>0\})$, the following three statements hold:
\begin{description}
\item (1) $F[x,y] < L^-[x,y]<L^+[x,y]$.
\item (2) For any $L>0$ in the open interval $ (L^-[x,y],L^+[x,y])$, there is not a feasible Dubins path with a length of $L$ from $\boldsymbol{z}_0$ to $(x,y)$.
\item (3) For any  $L>0$ in the closed interval $[F[x,y],L^-[x,y]]$ or in the semi-open interval $ [L^+[x,y],+\infty)$, there exists a feasible Dubins path with a length of $L$ from $\boldsymbol{z}_0$ to $(x,y)$.
\end{description}
\end{lemma}
This lemma is a direct result of \cite[Theorem 2]{Ding:2019}.  Thanks to Lemmas \ref{LE:continuity} and \ref{LE:Ding}, we have the following result.
\begin{theorem}\label{LE:R3}
If the minimum-time intercept point $E(t_m)$ between the pursuer and the target occurs in the interior of  $\mathcal{R}_3\cap\{x>0\}$, i.e., $E(t_m)\in \mathrm{int}(\mathcal{R}_3\cap \{x>0\})$, we then have
\begin{align}
\begin{split}
t_m =& \min \{ t>0 \big| t = F[E(t)],\  t ={L}^-[E(t)],\ \text{or}
  t={L}^+[E(t)]\}.
\end{split}
\label{EQ:tm_R3}
\end{align}
\end{theorem}
It can be seen from Theorem \ref{LE:R3} that the minimum intercept time $t_m$ may not be the fixed point of $F[E(t)]$. Up to present, it has not been clear what the solution type of the MTIP is if $t_m \neq F[E(t_m)]$.

By the following remark, we shall present the solution property of the MTIP for the cases of $t_m = L^-[E(t_m)]$ and $t_m = L^+[E(t_m)]$.
\begin{remark}\label{LE:type_L}
If the minimum-time intercept point $E(t_m)$ lies in $\mathrm{int}(\mathcal{R}_3)\cap \{x>0\}$,  the following statements hold:
\begin{description}
\item (1) If $t_m = L^-[E(t_m)]$, the solution of the MTIP is of type $L_{\alpha^-}R$.
\item (2) If $t_m = L^+[E(t_m)]$, the solution of the MTIP is of type $L_{\alpha^+}R$.
\end{description}
\end{remark}

%It should be noticed  from Theorem \ref{LE:R3} that  the solution of MTIP may not be the same as  that of RDP. This result was predicted in \cite{Mayer:2015} by presenting a counter example.  However, it is worth noting that even if the solution path of MTIP is not the same as that of RDP, the types of path keep being as CS or CC, as shall be presented in the next section.

%Not only do the above theorems present the necessary and sufficient conditions for the type  to be the same as that of RDP from $\boldsymbol{z}_0$ to $E(t_m)$, but also show the solution geometric properties of MTIP.

%In the case of $E(t_m)\in \mathrm{int}(\mathcal{R}_3)\cap \{x>0\}$, if $t_m \neq F[E(t_m)]$ so that the solution path of the MTIP is not the same as that of RDP from $\boldsymbol{z}_0$ to $E(t_m)$,

Note that the two types $L_{\alpha^-}R$ and $L_{\alpha^+}R$ belong to CC. Thus, combining Remark \ref{RE:R1_R2} with Remark \ref{LE:type_L}, we have the following result.
\begin{remark}
The solution of the MTIP lies in a sufficient family of 4 types in $\{RL,LR,LS,RS\}$.
\end{remark}

%As the solution of RDP is either CS or CC, it is concluded that the solution of the MTIP lies in a sufficient family of 4 types.

%Having the solution property of MTIP does not ensure to compute the solution of MTIP. In fact,

Theorems \ref{LE:Occurance_R1} and \ref{LE:R3} indicate that finding the minimum intercept time  $t_m$ is amount to computing the fixed points of $F[E(\cdot)]$, ${L}^-[E(\cdot)]$, and ${L}^+[E(\cdot)]$ over some specific intervals. Once $t_m$ is found, the final point $E(t_m)$ is available, and we can use the geometric properties in Remarks \ref{RE:R1_R2} and \ref{LE:type_L} to determine the solution path of the MTIP.  In the next section, an algorithm will be developed to find those fixed points so that the minimum interception time $t_m$ can be obtained efficiently.

%By extending \cite[Lemma 4]{Mayer:2015}, the following lemma holds without the assumption on the continuity of $F[E(t)]$.

%{\color{blue}****And also for the boundary point of $R_1$, it can be proven that $t_m \in  \{t> 0 | t = F[E(t)], L^-[E(t)], or  L^+[E(t)]\}$ **** }.

\section{Algorithm for the solution of MTIP}\label{SE:Algorithm}

In this section, we first present a robust and efficient algorithm to find the zeros of sufficiently smooth real-valued function in Subsection \ref{Subse:algorithm}, which will be employed in Subsection \ref{SE:procedure}  to establish numerical methods to find the solution of the MTIP for the case that the target's velocity $\boldsymbol{v}$ is constant.

For simplicity of presentation, the proofs for all the lemmas of this section are postponed to Appendix \ref{Appendix:B}.

\subsection{Algorithm for finding zeros of sufficiently smooth real-valued functions}\label{Subse:algorithm}

Before proceeding, we first present a lemma regarding the relationship between  extremas and zeros of a sufficiently smooth real-valued function.
\begin{lemma}\label{LE:extreme_to_zero}
Given a sufficiently smooth function $G(t)$ so that its number of zeros over an interval $[a,b]$ is finite, denote by $t_1$, $t_2,$, $\ldots$, $t_n$ in $[a,b]$ the zeros of the differentiation of $G(t)$ with respect to time, i.e., $G^{\prime}(t_i) = 0$ for $i=1,\ldots,n$. If $a\leq t_1 < t_2 < \ldots < t_n \leq b$, the following two statements hold:
\begin{description}
\item (1) if $G(t_i) \times G(t_{i+1}) >0$,  the function $G(t)$ on the interval $[t_i,t_{i+1}]$ does not have a zero;
\item (2) if $G(t_i) \times G(t_{i+1}) <0$, the function  $G(t)$  on the interval $[t_i,t_{i+1}]$   has  only one zero.
 \end{description}
\end{lemma}
Thanks to this lemma, if $G(t_i)\times G(t_{i+1})< 0$, we can use a simple bisection method to find the only zero in the interval $(t_i,t_{i+1})$. For notational simplicity, if $G(t_i)\times G(t_{i+1})< 0$, we denote by
$$z = \textbf{B}[G(t),t_{i},t_{i+1}]$$
the bisection method to find the zero $z$ of $G(t)$ in the interval $[t_i,t_{i+1}]$. Then, we can use Algorithm \ref{algorithm1} to compute all the zeros of $G(\cdot)$ in a constant time.

\begin{algorithm}
\caption{ }\label{algorithm1}
Given a sufficiently smooth function $G(t)$ so that its number of zeros is finite over $[a,b]$, let $t_1<t_2<\ldots<t_n$ be all the zeros of the differentiation of $G(t)$ with respect to $t$, i.e., $G^{\prime}(t_i) = 0$ for $i=1,2,\ldots,n$. Then, all the zeros of $G(t)$ over $[a,b]$ can be found by the following procedures:
\begin{description}
\item 1. set $i=0$, $t_0 = a$, $t_{n+1} = b$, and $Z=\varnothing$
\item 2. \textbf{while} $i<n$
\item 3. \qquad \textbf{if} $G(t_i) = 0$
%\begin{description}
%\item \qquad \qquad  $z_j = t_i$
\item 4. \qquad \qquad $Z = Z\cup \{t_i\}$
%\end{description}
\item 5. \qquad \textbf{elseif} $G(t_i)\times G(t_{i+1}) < 0$
%\begin{description}
\item 6. \qquad \qquad  $z =\textbf{B}[G(t),t_i,t_{i+1}]$
\item 7.  \qquad \qquad  $Z = Z\cup \{z\}$
%\end{description}
\item 8.  \qquad  \textbf{endif}
\item 9. \qquad  $i=i+1$
\item 10. \textbf{endwhile}
\end{description}
\end{algorithm}
Let us gather a few words to explain the pseudo codes in Algorithm \ref{algorithm1}. For any given $i\in \{0,1,\ldots,n+1\}$, if $G(t_i) = 0$, we have that $t_i$ is a zero of $G(t)$; thus, we add $t_i$ into the set $Z$ of zeros at step 4. If $G(t_i)\times G(t_{i+1}) < 0$, according to Lemma \ref{LE:extreme_to_zero} we have that the function $G(t)$ on the interval $(t_i,t_{i+1})$ has  only one zero. Thus, a typical bisection method can be used to find that zero, as shown by step 6, and the zero is added to the set $Z$ at step 7. If $G(t_i)\times G(t_{i+1}) > 0$, none zero exists between $t_i$ and $t_{i+1}$ according to the first statement of Lemma \ref{LE:extreme_to_zero}. Thus, nothing is done in the while loop if $G(t_i)\times G(t_{i+1}) > 0$. As a result, after the while loop,
the set $Z$ defined in Algorithm \ref{algorithm1} contains all the zeros of $G(t)$. It should be noted that the bisection method in step 6 can be completed within a constant time. To this end, given any sufficiently smooth real-valued function $G(t)$, if the number of zeros is finite over the interval $[a,b]$, we can use Algorithm \ref{algorithm1} to find all the zeros within a constant time.

%Although many algorithms are available to find zeros of a nonlinear equation, to the authors' best knowledge, there is not  an algorithm that can find all the zeros in the literature.

 In the following subsection, Algorithm \ref{algorithm1} will be applied to computing the solution of the MTIP.

%Let us gather a few words on how to apply Algorithm \ref{algorithm1} to find all the zeros of the formula in Eq.~\eqref{EQ:bar_Gcs} in the following remark.

%

\subsection{Computing the solution of the MTIP}\label{SE:procedure}

In this subsection, the geometric properties revealed in Section \ref{SE:Syntheses} will be employed to establish some nonlinear equations  so that  the solution length of the MTIP is determined by a specific zero of the nonlinear equations. As a result, Algorithm \ref{algorithm1} can be applied to finding the solution of the MTIP.

%The time for the target's path $E(t)$ to intersect with the boundaries of $\mathcal{R}_1$, $\mathcal{R}_2$, and $\mathcal{R}_3$ can be readily computed by simple geometric analysis. Therefore, this section focuses on presenting analytic methods to compute $t_m$ for the case that $E(t_m) $ lies in the interior of $\mathcal{R}_1$, $\mathcal{R}_2$, and $\mathcal{R}_3$.

Due to the symmetric property presented in Remark \ref{RE:Symmetric}, we will only consider the scenario that the final point $E(t_m)$ lies on the right plane $\{x>0\}$ in the following subsections.

\subsubsection{The case of $E(t_m) \in \mathrm{int}(\mathcal{R}_1)\cap\{x>0\}$} \label{subsection:R2}

According to Theorem \ref{LE:Occurance_R1}, if $E(t_m) \in \mathrm{int}(\mathrm{R}_1\cap\{x>0\})$, we have $t_m = F(t_m)$, indicating that the solution path is the same as that of RDP from $\boldsymbol{z}_0$ to $E(t_m)$ (cf. Lemma \ref{LE:fixed_RDP}). For any terminal point in $\mathcal{R}_1\cap \{x>0\}$, the solution path of RDP is of type RS \cite{Boissonnat:1994}, as shown by Fig.~\ref{Fig:Types_Region} in Appendix \ref{Appendix:A}. By the following lemma, a nonlinear equation in terms of the parameters of the RS path will be established, so that the solution of the MTIP can be found by finding a specific zero of the nonlinear equation.
\begin{lemma}\label{LE:analytic_CS}
Assume that the target's velocity is constant. If the minimum-time intercept point $E(t_m)$ between the pursuer and the target occurs in the interior of $\mathcal{R}_1\cap\{x>0\}$, i.e., $E(t_m)\in \mathrm{int}(\mathcal{R}_1\cap\{x>0\})$, we have that  the solution path of the MTIP is of type R$_{\alpha}$S and it holds
\begin{align}
G_{cs}(\alpha) \overset{\triangle}{=} \  &   {A_1 \sin \alpha + A_2 \cos \alpha}+ \alpha (A_3 \cos \alpha  + A_4  \sin \alpha ) + {A_5}= 0
  \label{EQ:CS_alpha}
\end{align}
where $\alpha\in [0,2\pi]$ is the radian of the right-turning circular arc R$_{\alpha}$, and  $A_1$--$A_5$ are constants in terms of $\rho$ and $\boldsymbol{z}_0$.
\end{lemma}

The expressions of $A_1$--$A_5$ are given in the proof of this lemma in Appendix \ref{Appendix:B}.

%Once the radian $\alpha\in [0,2\pi]$ of the right-turning circular arc R is avaialble, we are able to compute the solution by simple geometric analysis.
 Lemma \ref{LE:analytic_CS} shows that the radian $\alpha \in [0,2\pi]$ can be computed by finding the zeros of Eq.~\eqref{EQ:CS_alpha}. However,  Eq.~\eqref{EQ:CS_alpha} is a transcendental equation that may have multiple zeros.  Thus, the typical Newton-like iterative method or bisection method, proposed in \cite{McGee:2007}, may not find the desired zero that is related to the solution. In the following paragraph, a variant of $G_{cs}(\alpha)$ in Eq.~(\ref{EQ:CS_alpha}) will be presented so that Algorithm \ref{algorithm1} can be applied to finding the radian $\alpha$ of the right-turning circular arc.

%In the following paragraphs, some properties of sufficiently smooth functions will be presented so that a robust and efficient method will be developed to find all the zeros on the transcendental equation in Eq.~\eqref{EQ:CS_alpha}. To this end, one just needs to select the desired zero among all zeros in order to compute the solutoin of MTRP.

By rearranging  Eq.~\eqref{EQ:CS_alpha}, if $A_3 \cos \alpha +A_4 \sin \alpha \neq 0$ we have that $G_{cs}(\alpha)$ in Eq.~\eqref{EQ:CS_alpha}  is equivalent to
\begin{align}\label{EQ:bar_Gcs}
\bar{G}_{cs}(\alpha) \overset{\triangle}{ =}
\alpha +  \frac{A_1 \sin \alpha + A_2 \cos \alpha + A_5}{A_3 \cos \alpha + A_4 \sin \alpha}
\end{align}
Differentiating this equation with respect to $\alpha$ leads to
\begin{align}\label{EQ:dG_da}
\frac{\mathrm{d} \bar{G}_{cs}(\alpha)}{\mathrm{d} \alpha} =  1&+   \big[(A_1 \cos \alpha  - A_2 \sin \alpha )(A_3 \cos \alpha 
 +A_4 \sin \alpha) \notag\\
 &-   (A_1 \sin \alpha + A_2 \cos \alpha + A_5) 
 \times ( A_4 \cos \alpha - A_3 \sin \alpha )\big]/(A_3 \cos \alpha
  + A_4 \sin \alpha)^2
\end{align}
By substituting the half-angle formulas
\begin{align}\label{EQ:tan_x}
\sin \alpha = \frac{2 \tan \frac{\alpha}{2}}{ 1 + \tan^2 \frac{\alpha}{2}}\ \text{and}\ \cos \alpha = \frac{1 -  \tan^2 \frac{\alpha}{2}}{ 1 + \tan^2 \frac{\alpha}{2}}
\end{align}
into Eq.~\eqref{EQ:dG_da}, we have that  $\tan(\alpha/2)$  is a zero of the following  quartic polynomial:
\begin{align}
B_1 x^4 + B_2 x^3 + B_3 x^2 + B_4 x + B_5 = 0
\label{EQ:polynomial_quartic}
\end{align}
where
\begin{align}
\begin{split}
B_1 & =-A_3+A_1A_3+A_4A_5-A_2A_4 \\
B_2 & = 2A_4+2A_3A_5\\
B_3 &=2A_1A_3-3A_2A_4+A_4A_5+A_2-A_5\\
B_4 & = 2A_4+2A_3A_5\\
B_5 & =A_3+A_1A_3-A_2A_4-A_4A_5
\end{split}\nonumber
\end{align}
The roots of any quartic polynomial can be readily obtained either by radicals or by standard polynomial solvers.
% Thus, the zeros of the formula in Eq.~\eqref{EQ:dG_da} can be found efficiently.
By the following remark, we shall show how to apply Algorithm \ref{algorithm1} to finding the shortest path of type R$_{\alpha}$S from the zeros of the quartic polynomial in Eq.~(\ref{EQ:polynomial_quartic}).
%By the following lemma, we shall show that the zeros of the original equation in Eq.~\eqref{EQ:bar_Gcs} can be efficiently found from the zeros of Eq.~\eqref{EQ:dG_da}.

\begin{remark}\label{RE:procedure}
If $E(t_m)\in \mathrm{int}(\mathcal{R}_1\cap \{x>0\})$ so that the solution path of the MTIP is of type $R_{\alpha}S$, then we can use the following procedure to find all the zeros of Eq.~\eqref{EQ:bar_Gcs}:
\begin{description}
\item (1) find all the real zeros of the quartic polynomial in Eq.~\eqref{EQ:polynomial_quartic} by either radicals or by a standard polynomial solver;
\item (2) by combining  Eq.~\eqref{EQ:tan_x} and the real zeros of Eq.~\eqref{EQ:polynomial_quartic}, we can find all the real zeros of   Eq.~\eqref{EQ:dG_da};
\item (3) because all the real zeros of Eq.~\eqref{EQ:dG_da} are the extremas of Eq.~\eqref{EQ:bar_Gcs} and the function $\bar{G}_{cs}(\alpha)$ satisfies the conditions of Lemma \ref{LE:extreme_to_zero}, it follows that Algorithm \ref{algorithm1} can be used to find all the zeros of Eq.~\eqref{EQ:bar_Gcs} efficiently;
\item (4) if we denote by $\alpha_1$, $\alpha_2$, $\ldots$, $\alpha_m$ the zeros of Eq.~\eqref{EQ:bar_Gcs}, any path of type R$_{\alpha_i}$S for $i=1,2,\ldots,m$ can be computed by a simple geometric analysis, and the shortest path of the type R$_{\alpha_i}$S is the solution of the MTIP.
\end{description}
\end{remark}

%If we consider each zero of Eq.~\eqref{EQ:bar_Gcs} as the radian of the right turning circular arc, we are able to readily obtain a path of type R$_{\alpha}$S by simple geometric analysis. As a result, we

%Compute the length of solution path related to each computed $\alpha$ by simple geometric analysis, and select the specific $\alpha$   related to the shortest path as the desired radian of the right-turning circular arc.

%If $E(t_m)\in \mathrm{int}(\mathcal{R}_1\cap \{x<0\})$, we can combine the symmetric property in Remark \ref{RE:Symmetric} and the procedure in Remark \ref{RE:procedure} to find the solution of the MTIP.

%Once the radian $\alpha$ of the right-turning circular arc is obtained, we are able to compute the solution path of type $R_{\alpha}S$ immediately from geometric analysis.

\subsubsection{The case of $E(t_m) \in \mathrm{int}(\mathcal{R}_2\cap\{x>0\})$}\label{Subsection:R1}

In view of Theorem \ref{LE:Occurance_R1}, if $E(t_m) \in \mathrm{int}(\mathrm{R}_2)$, we have $t_m = F(t_m)$, indicating that the solution path is the same as that of the RDP from $\boldsymbol{z}_0$ to $E(t_m)$ (cf. Lemma \ref{LE:fixed_RDP}). For any terminal point in $\mathcal{R}_2\cap\{x>0\}$, the solution path of the RDP is of type LR \cite{Boissonnat:1994}, as shown by Fig.~\ref{Fig:Types_Region} in Appendix \ref{Appendix:A}.

 By the following lemma, an equation in terms of the parameters of LR will be established.
\begin{lemma}\label{LE:analytic_CC}
If the minimum-time intercept point $E(t_m)$ occurs in $\mathrm{int}(\mathcal{R}_2\cap\{x>0\})$,  then the solution of the MTIP is of type $L_{\beta}R_{\gamma}$ and it holds
\begin{align}\label{EQ:G_cc}
G_{cc}(\eta) \overset{\triangle}{=} \ & B_1 \eta^4 + B_2 \eta^3 + B_3 \eta^2 + B_4 \eta + B_5 +  B_6 \cos \eta  +  B_7 \sin \eta + \eta (B_8 \cos \eta + B_9 \sin \eta)  = 0
\end{align}
where $\eta = \beta + \gamma$ is  the sum of the radians of the left-turning circular arc $L_{\beta}$ and the right-turning circular arc $R_{\gamma}$, and  $B_1$--$B_9$ are constants in terms of $\rho$ and $\boldsymbol{z}_0$.
\end{lemma}
The expressions of $B_1$--$B_9$ are given in the proof of this lemma in Appendix \ref{Appendix:B}.

The fourth derivative of $G_{cc}$ is expressed as
\begin{align}\label{EQ:G_cc4}
G_{cc}^{(4)}(\eta) =&\ (B_7 + 4B_8) \sin \eta + (B_6 - 4 B_9) \cos \eta +\eta (B_8 \cos \eta + B_9 \sin \eta) + 24 B_1.
\end{align}
Note  that the form of $G_{cc}^{(4)}(\eta)$ is the same as that of $G_{cs}(\alpha)$ in Eq.~\eqref{EQ:CS_alpha}. Thus, all the zeros of  $G_{cc}^{(4)}(\eta)$ can be obtained in a constant time by the procedure in Remark \ref{RE:procedure}.  According to Lemma \ref{LE:extreme_to_zero},  Algorithm \ref{algorithm1} can be applied to finding the zeros of $G^{(i)}_{cc}(\eta)$ from the zeros of $G^{(i+1)}_{cc}(\eta)$. Hence, we are able to find all the real zeros of $G_{cc}(\eta)$ in Eq.~\eqref{EQ:G_cc} by applying Algorithm \ref{algorithm1} four times. Note that we have $t_m = \rho \eta$. Thus, the minimum positive zero of Eq.~\eqref{EQ:G_cc} gives the minimum interception time $t_m$.
%one can immediately select the one related to the solution of MTIP. Finally, the minimum time for interception or the length of the solution path of MTIP is given by $\rho \eta$, i.e., $t_m = \rho \eta$.

\subsubsection{The case of $E(t_m) \in \mathrm{int}(\mathcal{R}_3)\cap \{x> 0\}$}\label{Subse:R3}

If $E(t_m) \in \mathrm{int}(\mathcal{R}_3\cap \{x>0\})$, not only do we need to check the fixed points of $F[E(t)]$ but also the fixed points of $L^-[E(t)]$ and $L^+[E(t)]$ according to Theorem \ref{LE:R3}. Since the path related to $F[E(t)]$ in $\mathrm{int}(\mathrm{R}_3\cap\{x>0\})$ is of type $RS$ \cite{Boissonnat:1994}, it follows that  the fixed points of $F[E(t)]$ in this case can be obtained by the procedure in Remark \ref{RE:procedure}.

Even though the paths related to ${L}^-[E(t)]$ and ${L}^+[E(t)]$ are not the same as the solution of  RDP from $\boldsymbol{z}_0$ to $E(t_m)$, their types are LR (cf. Lemma \ref{LE:type_L}). If we denote by $\eta> 0$ the sum of the radians of L and R, we have that Eq.~\eqref{EQ:G_cc} holds for both $t_m = L^+[E(t_m)]$ and $t_m=L^-[E(t_m)]$ according to the proof of Lemma \ref{LE:analytic_CC}. Therefore, the fixed points of ${L}^-[E(t)]$ and ${L}^+[E(t)]$ can be obtained by the same procedure as presented in Subsection \ref{Subsection:R1}.

%the proof of Lemma that the fixed points of $\mathcal{L}_1(\cdot)$ and $\mathcal{L}_2(\cdot)$ are the zeros of $G_{cc}$.

%From Lemma \ref{lemma4}, we know that there are four types of Relaxed-Dubins paths and the optimal one has the shortest length. So in this section, we derive all kinds of semi-analytical solutions, including LS, RS, LR, RL, as follows.

%\subsection{Semi-analytical Solution for LS}

\section{Numerical examples}\label{SE:Examples}

In this section, we present some examples to demonstrate the developments of the paper. In Subsection \ref{Subse:MTIP}, four examples of the MTIP are simulated, and in Subsection \ref{Subse:drift} an example is present to show how to apply the above algorithm to planning Dubins path in a constant drift field.

Note that the position vector $(x,y)$ is normalized so that the speed of pursuer is one in Section \ref{SE:Preliminary}. Thus, the units of variables in this section are omitted.

\subsection{Examples of the MTIP}\label{Subse:MTIP}

\subsubsection{Case A}

\begin{figure}[h]
	\centering
	\includegraphics[width = 9cm]{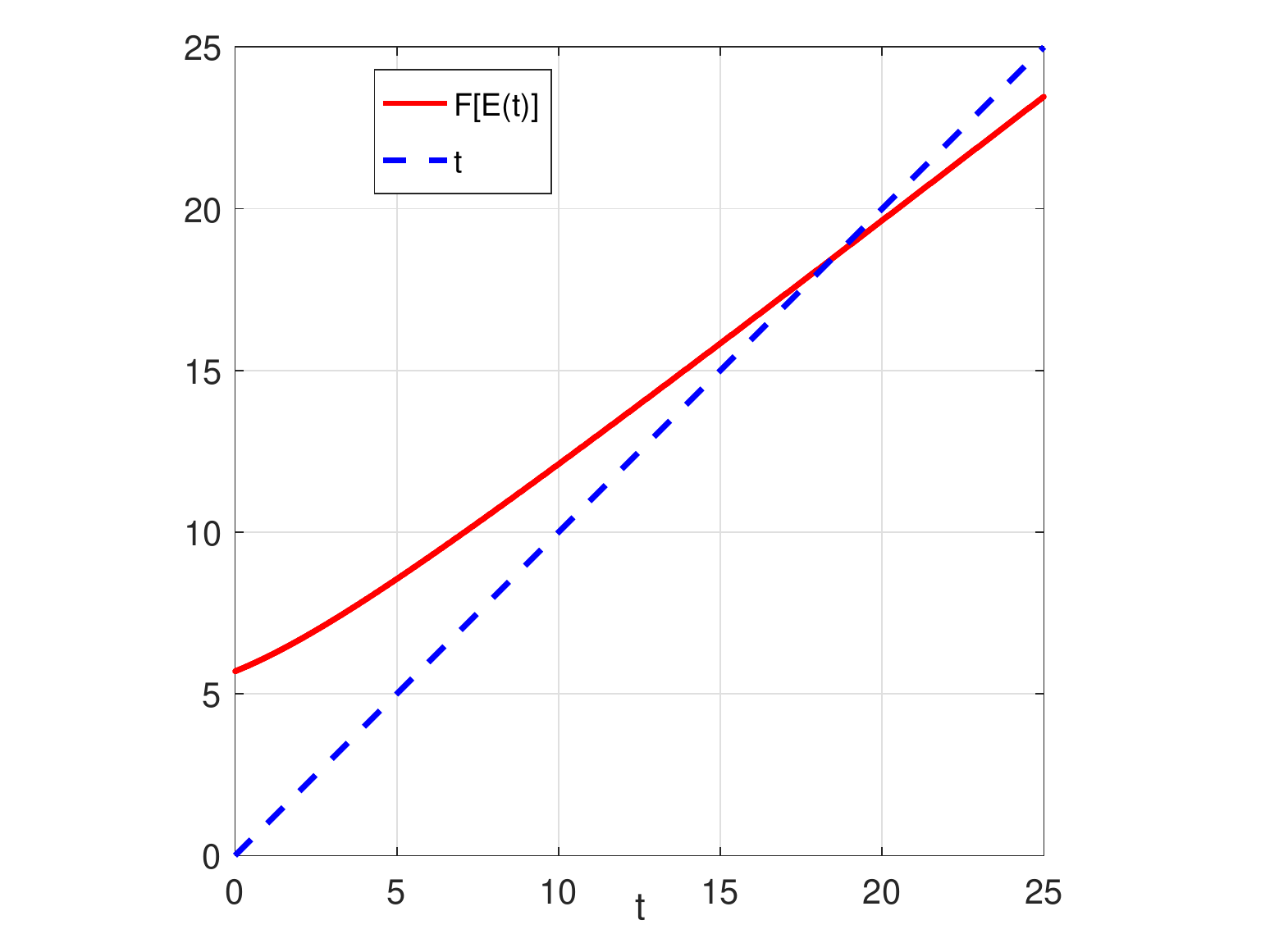}
	\caption{Case A: The profile of $F[E(t)]$ with respect to time.}
	\label{Fig:CaseAA1}
\end{figure}

For case A, we present an example for which the function $F[E(t)]$ is continuous. We set
$\rho = 1$, $(\hat{x}_{0},\hat{y}_{0}) = (5,2)$, and $\boldsymbol{v} = (0.55,-0.55)$ to ensure that the half line $E(\cdot):[0,+\infty]$ does not intersect the two half circles $\{(x,y)\in C_r \cup C_l\ |\ y>0 \}$. According to Lemma \ref{TH:continuity}, these parameters guarantee that the function $F[E(t)]$ is continuous, as shown in Fig.~\ref{Fig:CaseAA1}.

\begin{figure}[h]
	\centering
	\includegraphics[width = 9cm]{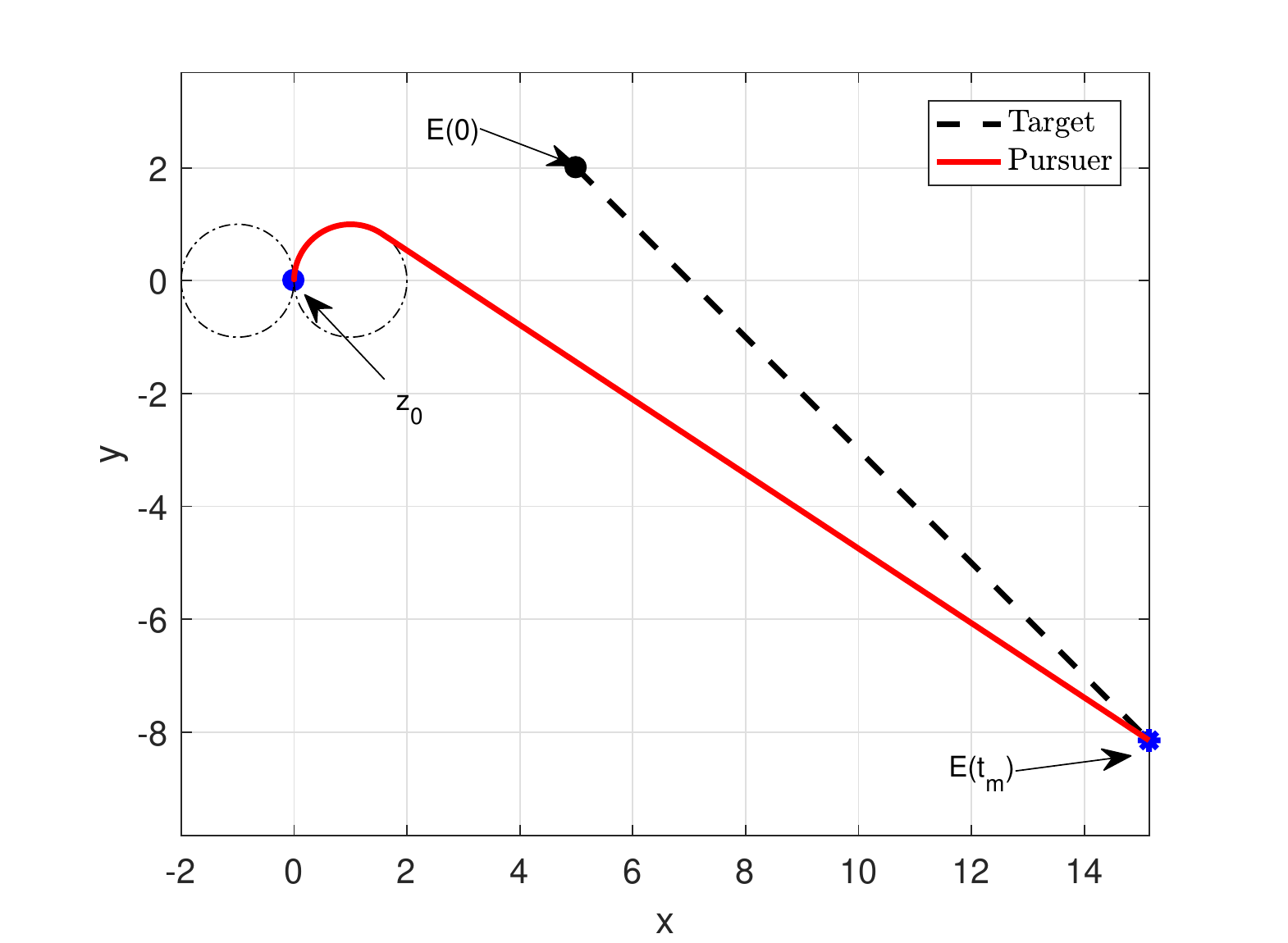}
	\caption{Case A: The solution path of the MTIP.}
	\label{Fig:CaseAA2}
\end{figure}

In such a case, the minimum intercept time $t_m$ is the minimum fixed point of $F[E(t)]$ according to Lemma \ref{LE:continuous}. By the developments in Subsections \ref{subsection:R2} and \ref{Subsection:R1}, all the fixed points of $F[E(t)]$ can be computed in a constant time. The time to compute the solution of case A is tested by MATLAB on a desktop with Intel(R) Core(TM)i3-4130U CPU@0.725GHz, showing that the solution is computed within $10^{-4}$ seconds. The minimum fixed point is computed as $18.45$, indicating that the minimum intercept time is $t_m = 18.45$. Since the target's velocity is constant, we have $E(t_m) = (\hat{x}_{0},\hat{y}_{0}) + \boldsymbol{v} t_m = (15.15,-8.15)$. In view of Lemma \ref{LE:fixed_RDP}, the solution path of the MTIP is the same as that of the RDP from $\boldsymbol{z}_0$ to $E(t_m) = ((15.15,-8.15))$. Thus, the solution of the MTIP for case A is readily available by geometric analysis \cite{Boissonnat:1994}, and the solution path is shown in  Fig.~\ref{Fig:CaseAA2}.

%By a simple coordinate transformation, the MTIP in case A is equivalent to the RDP in a constant drift field. The solution in the constant drift field is presented in Fig.~\ref{Fig:case1b}.
%\begin{figure}[h]
%	\centering
%	\includegraphics[width = 8cm]{case1b.eps}
%	\caption{The minimum interception path in the equivalent constant drift field.}
%	\label{Fig:case1b}
%\end{figure}

\subsubsection{Case B}

\begin{figure}[h]
	\centering
	\includegraphics[width = 9cm]{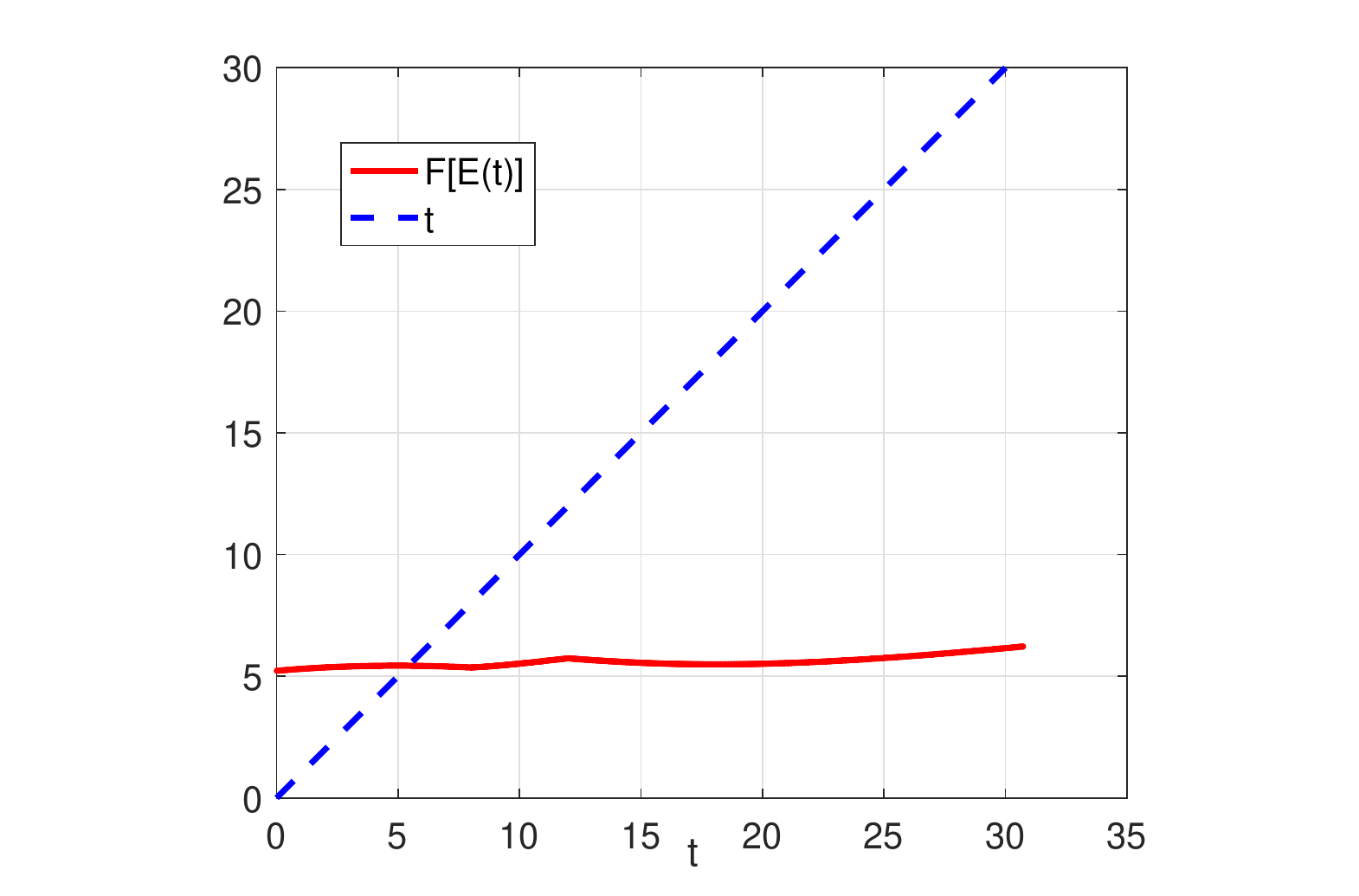}
	\caption{Case B: The profile of $F[E(t)]$ with respect to time.}
	\label{Fig:CaseBB1}
\end{figure}

\begin{figure}[h]
	\centering
	\includegraphics[width = 9cm]{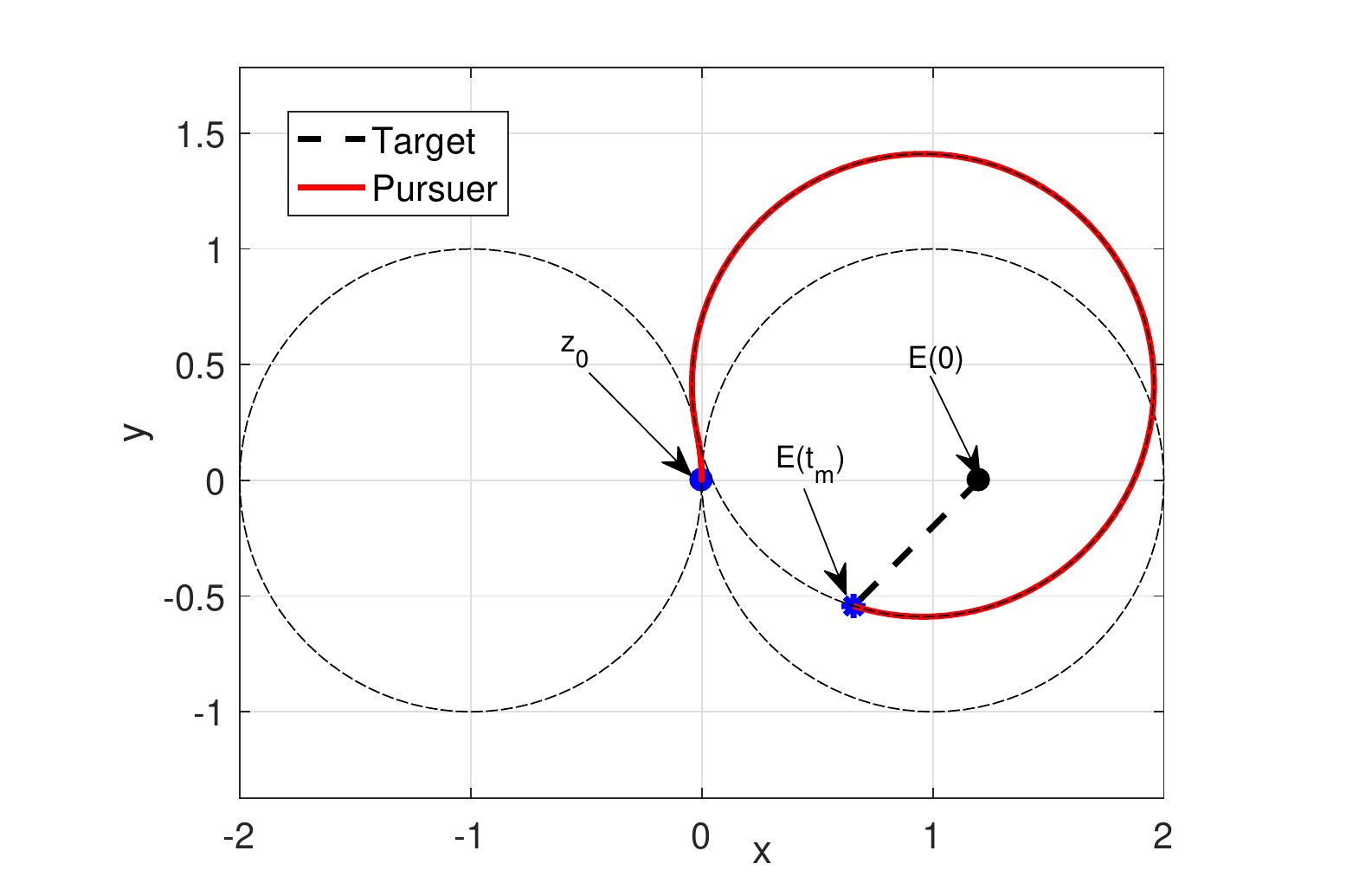}
	\caption{Case B: The solution path of the MTIP.}
	\label{Fig:CaseBB2}
\end{figure}

The initial conditions of case A are set so that the minimum-time interception point $E(t_m)$ between the pursuer and the target occurs in $\mathcal{R}_1$. For case B, we
choose $\rho = 1$, $(\hat{x}_{0},\hat{y}_{0})  = (1.2,0)$, and $\boldsymbol{v} = (-0.1,-0.1)$. These initial parameters are designed so that the minimum-time interception happens in $\mathcal{R}_2$. In addition, these initial parameters ensure that the half line $E(\cdot):[0,+\infty]$ does not intersect the two half circles $\{(x,y)\in C_r \cup C_l\ |\ y>0 \}$. Thus, according to Lemma \ref{TH:continuity}, the function $F[E(t)]$ for $t\geq 0$ is continuous, and it is plotted in Fig.~\ref{Fig:CaseBB1}
According to Lemma \ref{LE:continuous}, the minimum fixed point of $F[E(t)]$ is the minimum time for the interception between the pursuer and the target, and it is computed as $5.43$. Then, the interception point is $E(t_m) = (\hat{x}_{0},\hat{y}_{0}) + \boldsymbol{v}  t_m = (0.66,-0.54)$.  With this terminal point, the solution path of the MTIP is the same as that of the RDP from $\boldsymbol{z}_0$ to $E(t_m)$, and the solution path is readily available according to \cite{Boissonnat:1994}, as shown in Fig.~\ref{Fig:CaseBB2}.
The time to compute the solution of case B is tested by MATLAB on a desktop with Intel(R) Core(TM)i3-4130U CPU@0.725GHz, showing that the solution is computed within $10^{-4}$ seconds.

 %The path of minimum time interception is shown in Fig.~\ref{Fig:case4a}. For this case, the half line $E(\cdot):[0,t_f]$ intersects $R_1$ and $R_2$, so there is no ${L}^-[E(t)]$ and ${L}^+[E(t)]$ in Fig.~\ref{Fig:case4c} and the function $F(t)$ is continuous, so the only intersect point in Fig.~\ref{Fig:case4c} is the minimum-time interception point. It is also clear from Fig.~\ref{Fig:case4b} that the path for minimum-time interception in equivalent constant drift field.

%\begin{figure}[h]
%	\centering
%	\includegraphics[width = 8cm]{case4b.eps}
%	\caption{The minimum interception path in the equivalent constant drift field.}
%	\label{Fig:case4b}
%\end{figure}

\subsubsection{Case C}

\begin{figure}[h]
\centering
\includegraphics[trim = 9cm 0cm 6cm 0cm, clip, width = 9cm]{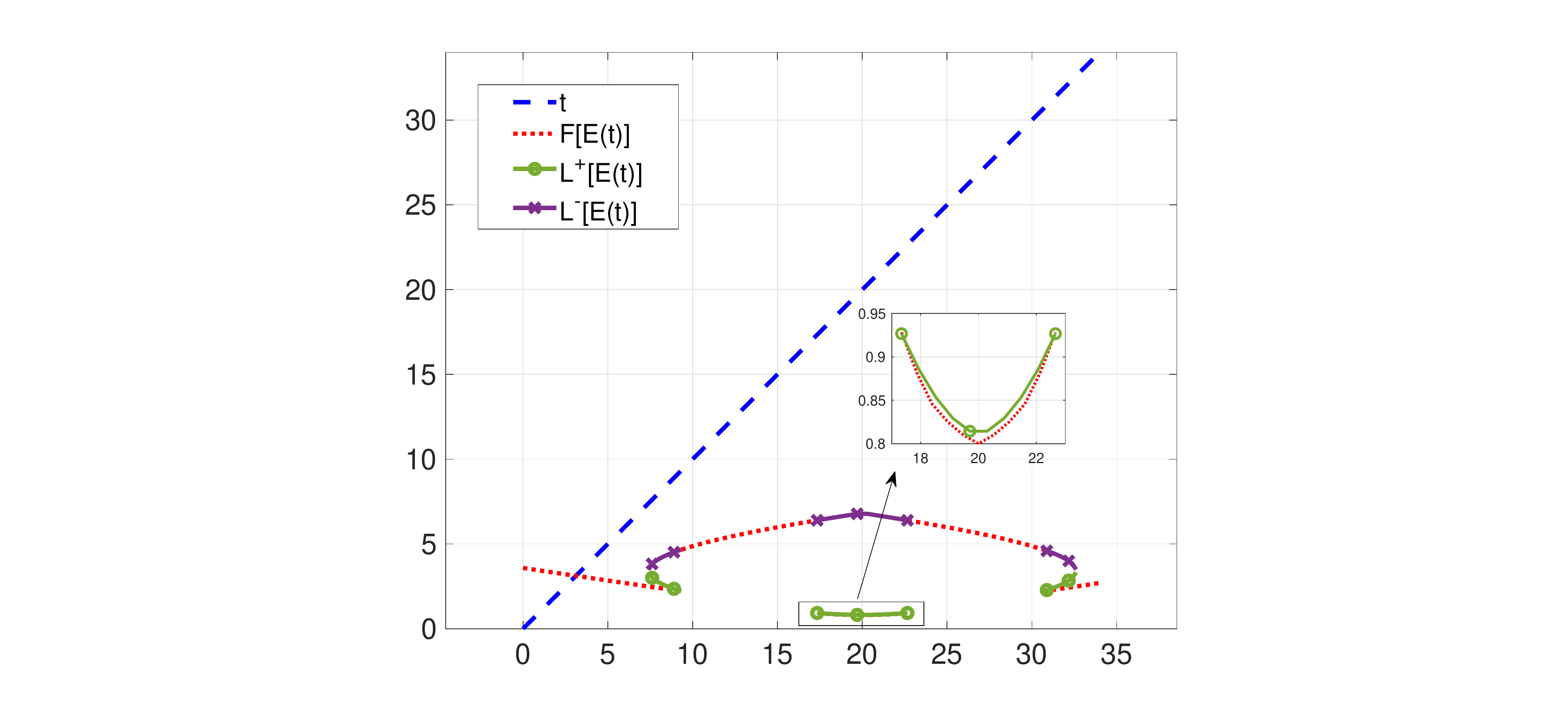}
\caption{Case C: The profiles of $F[E(t)]$, $L^-[E(t)]$, and $L^+[E(t)]$ with respect to time.}
	\label{Fig:CaseC1}
\end{figure}

For case C, we set $\rho = 1$, $(\hat{x}_{0},\hat{y}_{0})  = (-3,0.8)$, and $\boldsymbol{v} = (0.15,0)$. These initial parameters are tailored so that the half line $E(\cdot):[0,t_f]$ intersects the two half circles $\{(x,y)\in C_r \cup C_l\ |\ y>0 \}$. In such a case, the function $F[E(t)]$ is not continuous (cf. Lemma \ref{TH:continuity}), as plotted in Fig.~\ref{Fig:CaseC1}.
It is seen from Fig.~\ref{Fig:CaseC1} that the functions $L^-[E(t)]$ and $L^+[E(t)]$ exist over some intervals. This is reasonable because the target's path passes the subregion $\mathcal{R}_3$ over different intervals.

In the case that the function $F[E(t)]$ is discontinuous, checking the fixed points of $F[E(t)]$ is not enough to find the solution of the MTIP. One should also find the fixed points of $L^-[E(t)]$ and $L^+[E(t)]$. The procedures in Subsection \ref{SE:procedure} are applied, showing that the two functions $L^-[E(t)]$ and $L^+[E(t)]$ do not have a fixed point over their domains of definition. And, there is only one fixed point for the function $F[E(t)]$, which is computed as 3.15, i.e., $t_m = 3.15$.

\begin{figure}[h]
	\centering
	\includegraphics[width = 9cm]{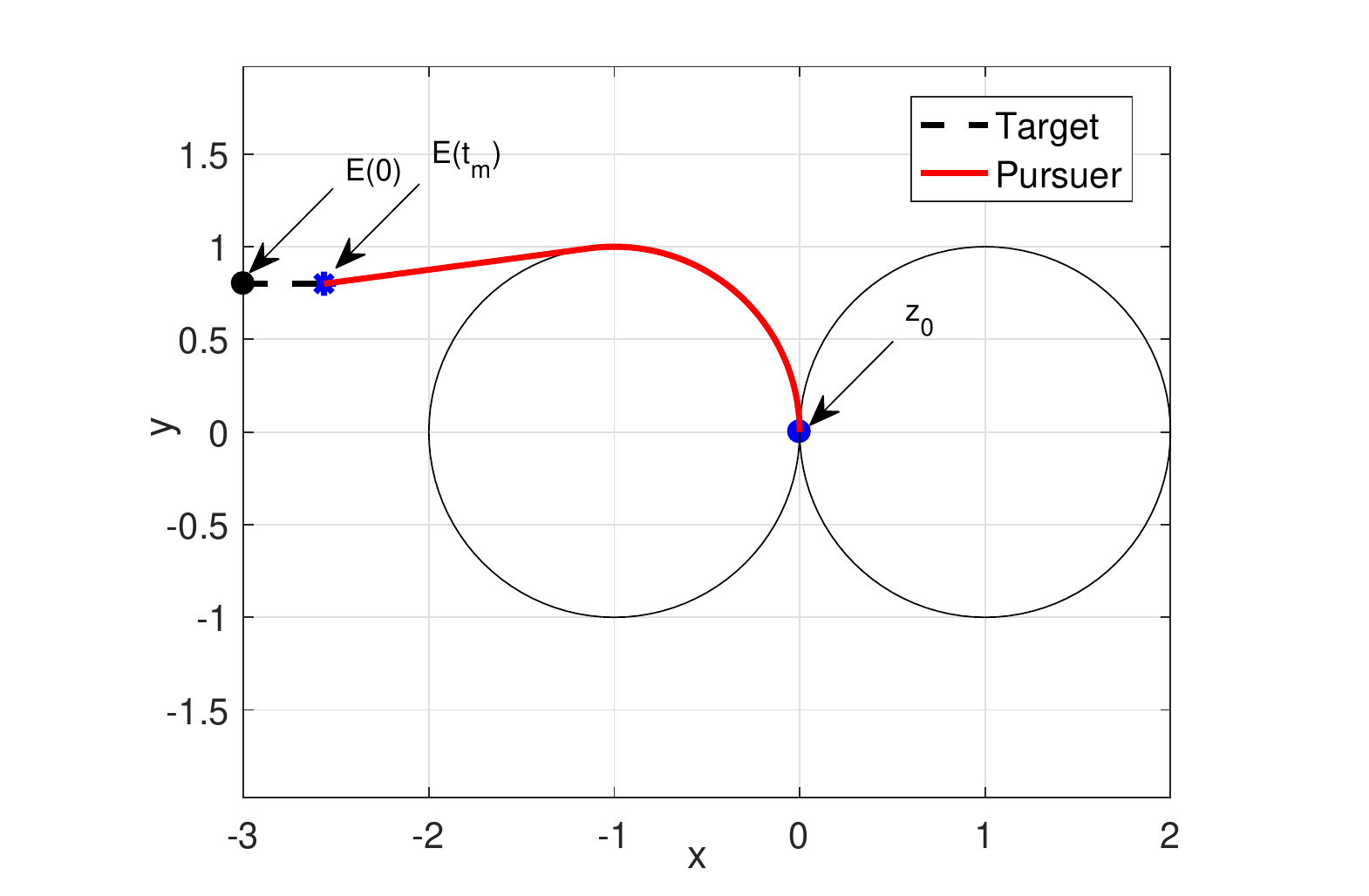}
	\caption{Case C: The solution path of the MTIP.}
	\label{Fig:CaseC2}
\end{figure}

Then, the interception point is $E(t_m) = (\hat{x}_{0},\hat{y}_{0}) + \boldsymbol{v} t_m = (-2.55,0.80)$.  With this terminal point, the solution path of the MTIP  is the same as that of the RDP from $\boldsymbol{z}_0$ to $E(t_m)$ according to Lemma \ref{LE:fixed_RDP}.The solution path is of type LS, as shown in Fig.~\ref{Fig:CaseC2}.
The time to compute the solution of case C is tested by MATLAB on a desktop with Intel(R) Core(TM)i3-4130U CPU@0.725GHz, showing that the solution is computed within $10^{-4}$ seconds.

%Fig.~\ref{Fig:case2a} shows the path of the minimum-time interception, and the path for minimum-time interception in equivalent constant drift field is shown in Fig.~\ref{Fig:case2b}.

%\begin{figure}[h]
%	\centering
%	\includegraphics[width = 8cm]{case2b.eps}
%	\caption{The minimum interception path in the equivalent constant drift field.}
%	\label{Fig:case2b}
%\end{figure}

\subsubsection{Case D}

All the solution paths in the above three cases are the same as those of RDP from $\boldsymbol{z}_0$ to $E(t_m)$. In \cite{Mayer:2015}, an example was presented to show that the solution of the MTIP was not the same as that of the RDP from $\boldsymbol{z}_0$ to $E(t_m)$. However, it was not shown how to find the solution path. Thanks to
the developments in Sections \ref{SE:Syntheses} and \ref{SE:Algorithm}, the example can be addressed efficiently as shown by the following paragraph.

\begin{figure}[h]
	\centering
	\includegraphics[width = 9cm]{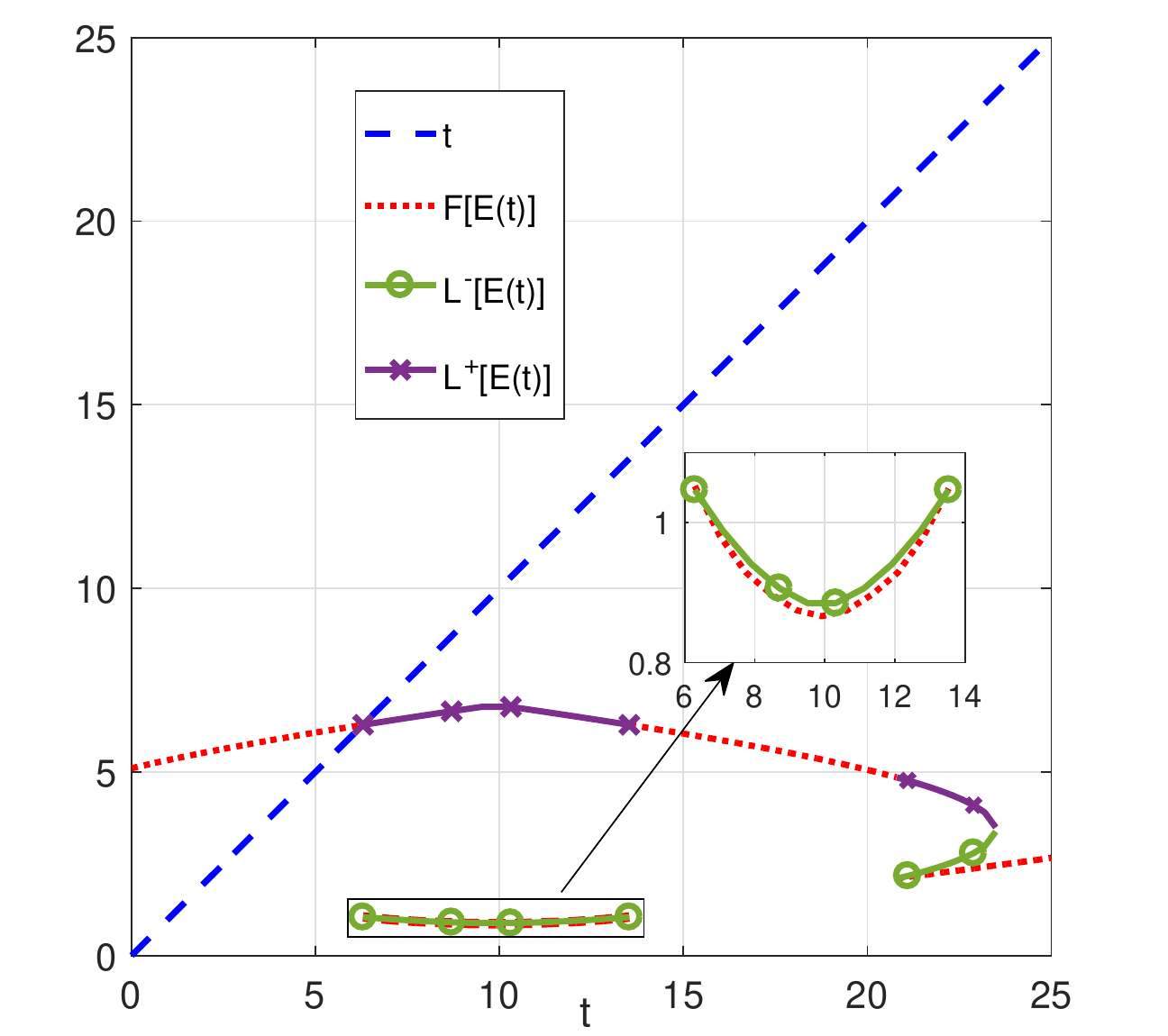}
	\caption{Case D: The profiles of $F[E(t)]$, $L^-[E(t)]$, and $L^+[E(t)]$ with respect to time $t$.}
	\label{Fig:CaseA2}
\end{figure}

For the example in \cite{Mayer:2015}, the minimum-turning radius is set as $\rho = 1$. The velocity vector and initial position of the target are given by $\boldsymbol{v} = (\frac{\sqrt{3}}{4\pi},0)$ and $(\hat{x}_{0},\hat{y}_{0})  = (-\frac{\sqrt{3}+1}{2},\frac{\sqrt{3}}{2})$, respectively.
 The profiles of $F[E(t)]$, $L^-[E(t)]$, and $L^+[E(t)]$ are presented in Fig.~\ref{Fig:CaseA2}.
It is seen from Fig.~\ref{Fig:CaseA2} that the function $F[E(t)]$ is discontinuous. The procedures in Subsection \ref{SE:procedure} are directly applied to computing the solution of the MTIP, showing that the minimum intercept time $t_m$ is a fixed point of $L^+[E(t)]$. Then, according to Remark \ref{LE:type_L}, the solution path is of type L$_{\alpha^+}$R, which is consistent with the analysis in \cite{Mayer:2015}.

%This is consistent with the analysis in \cite{Mayer:2015}.

 %This means that the final point $E(t_m)$ lies on the left turning circle $\mathcal{C}_l$.

\begin{figure}[h]
	\centering
	\includegraphics[width = 9cm]{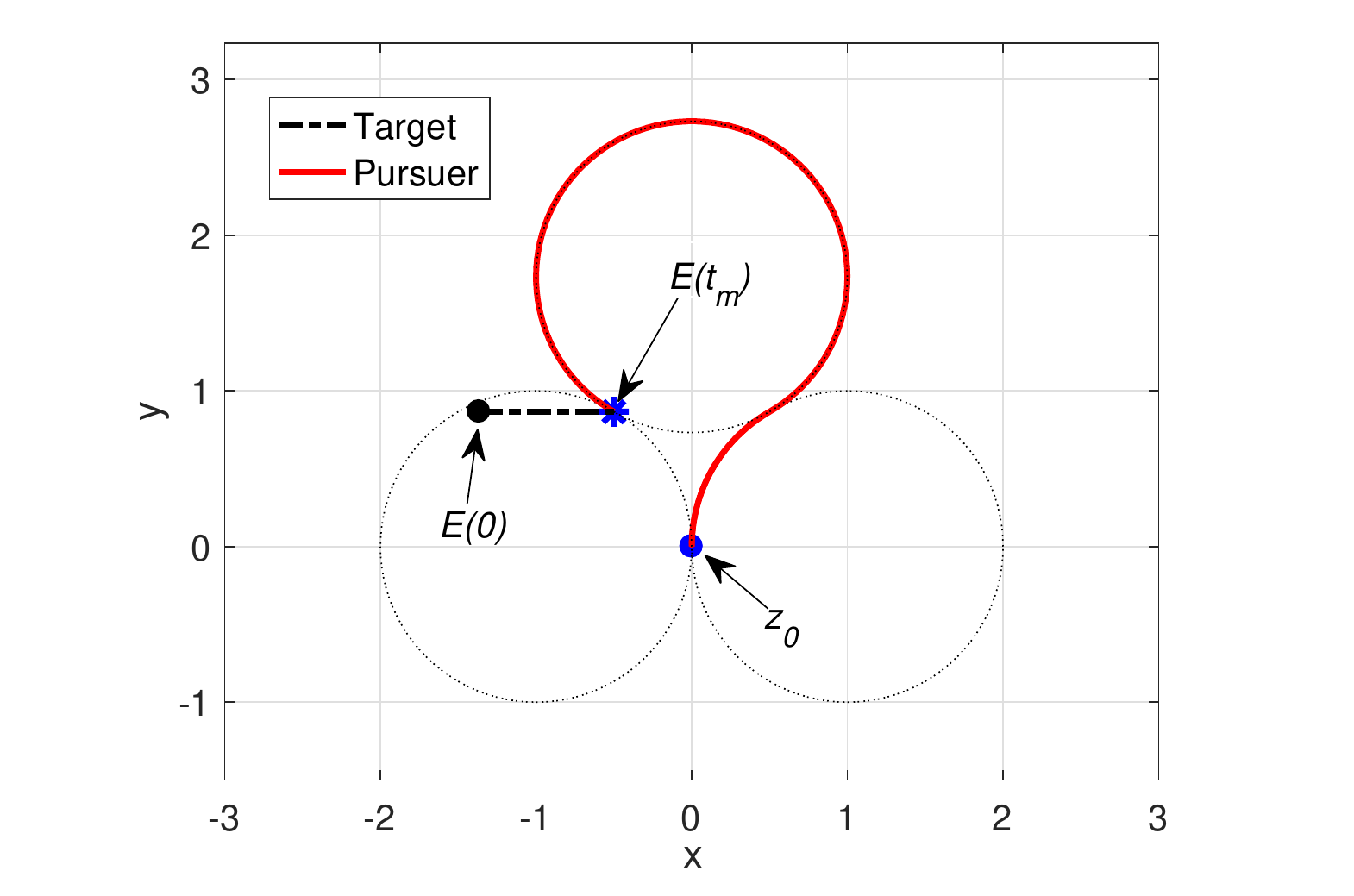}
	\caption{Case D: the solution path of the MTIP.}
	\label{Fig:CaseA1}
\end{figure}

Having the value of $t_m$, the solution path of the MTIP for this example can be computed by applying the results in Subsection \ref{Subse:R3}, and it is  presented in Fig.~\ref{Fig:CaseA1}.
It is apparent that the solution of RDP from $\boldsymbol{z}_0$ to $E(t_m)$ should be a single left turning circular arc, but the solution of the MTIP for case D is quite different. The time to compute the solution of case D is tested by MATLAB on a desktop with Intel(R) Core(TM)i3-4130U CPU@0.725GHz, showing that the solution is computed within $10^{-4}$ seconds.

\subsection{Path planning in constant drift field}\label{Subse:drift}

\begin{figure}[!htp]
\centering
\begin{subfigure}[t]{6cm}
\centering
\includegraphics[width = 6cm]{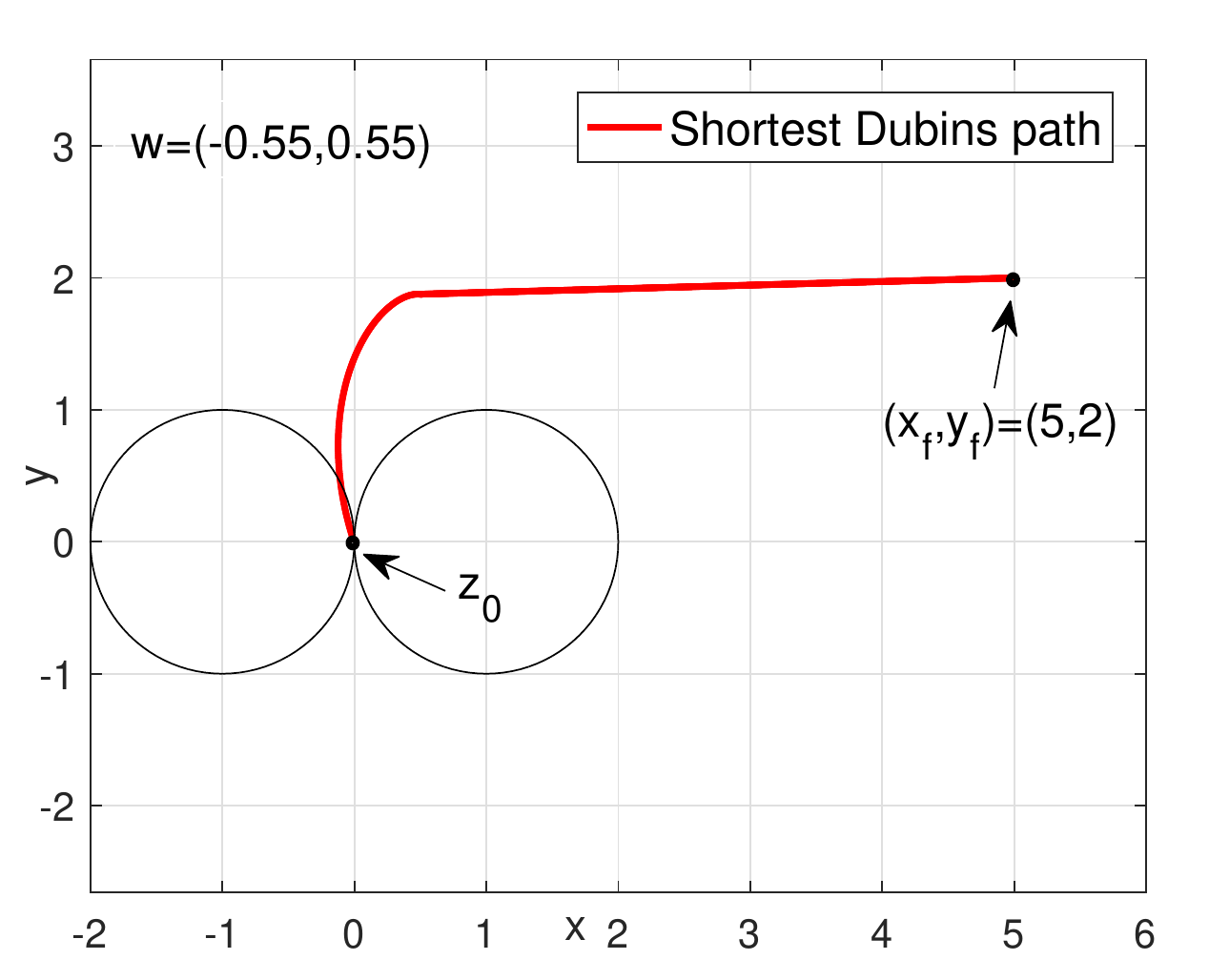}
%$\boldsymbol{w} = (-0.55,0.55)$\\
%$(x_f,y_f) = (5,2)$
\caption{}
\label{Fig:caseE1}
\end{subfigure}
\begin{subfigure}[t]{6cm}
\centering
\includegraphics[width = 6cm]{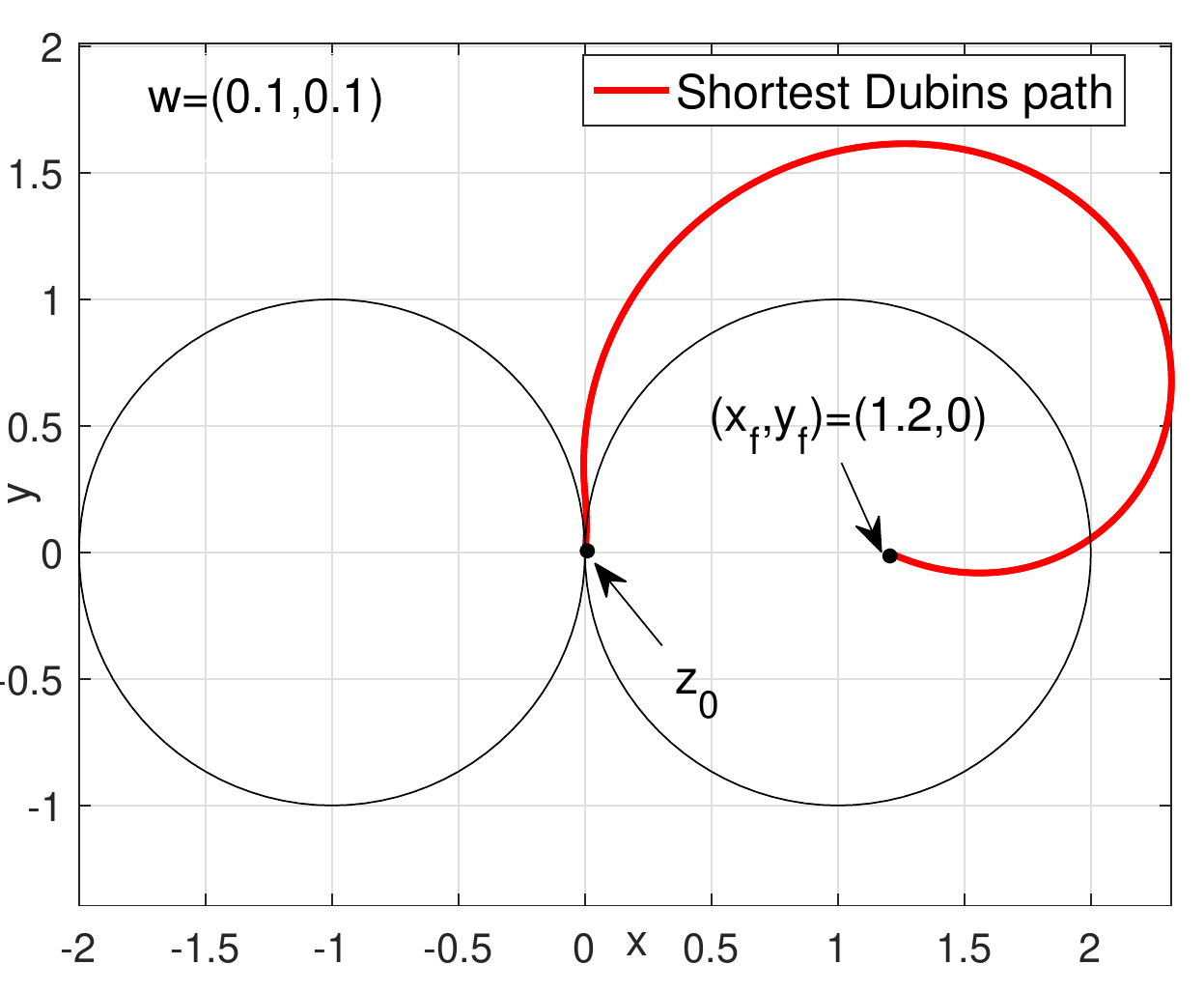}
%$\boldsymbol{w} = ()$ \\
% $(x_f,y_f) = ()$
\caption{}
\label{Fig:caseE2}
\end{subfigure}

\begin{subfigure}[t]{6cm}
\centering
\includegraphics[width = 6cm]{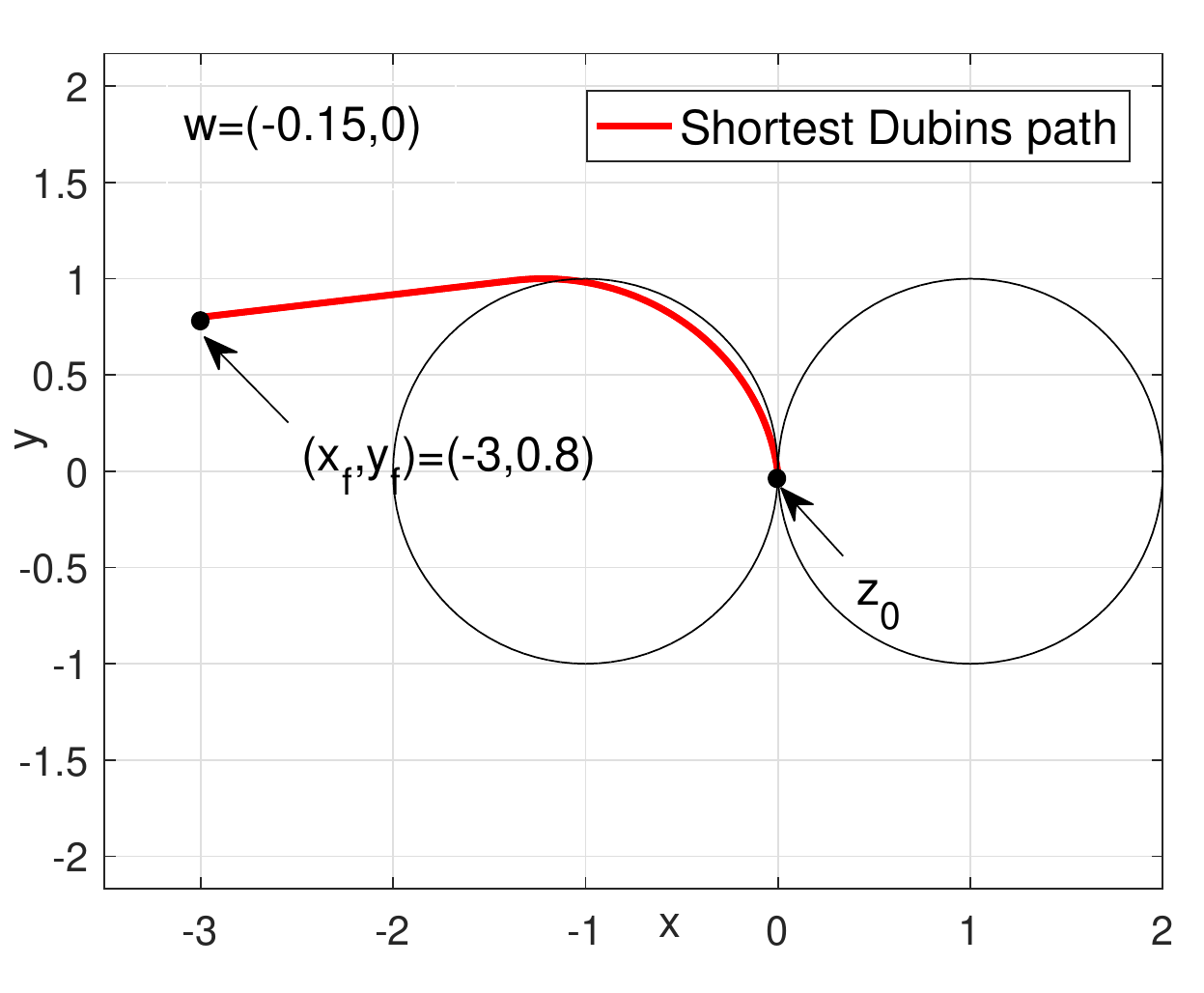}
%$\boldsymbol{w} = ()$ \\
% $(x_f,y_f) = ()$
\caption{}
\label{Fig:caseE3}
\end{subfigure}
\begin{subfigure}[t]{6cm}
\centering
\includegraphics[width = 6cm]{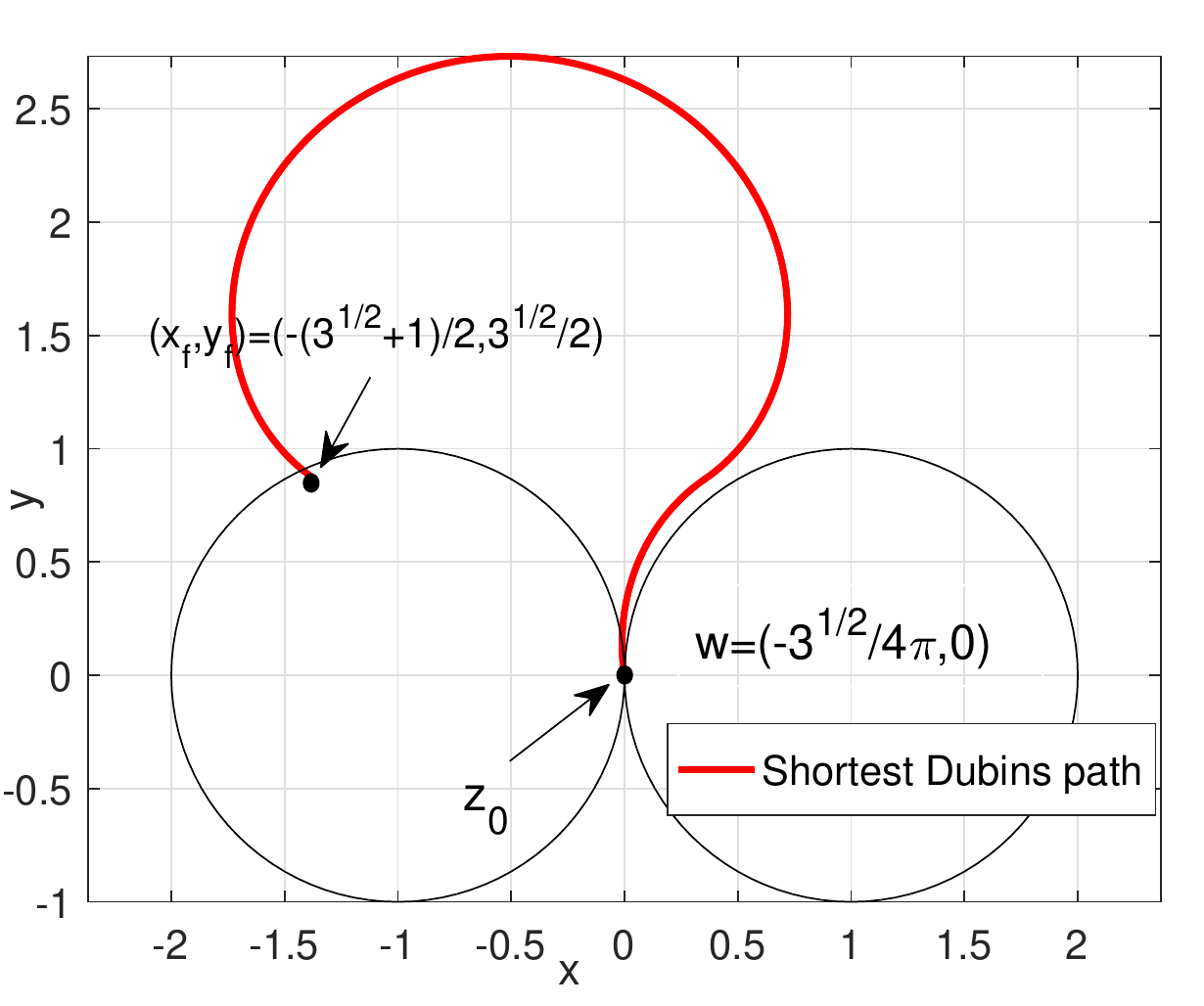}
%$\boldsymbol{w} = ()$ \\
% $(x_f,y_f) = ()$
\caption{}
\label{Fig:caseE4}
\end{subfigure}
\caption{Shortest Dubins paths in constant drift fields.}
\label{Fig:CaseE}
\end{figure}

Consider an aerial/marine vehicle whose motion is described by
\begin{align}
\dot{\boldsymbol{z}}(t) =
\left[
\begin{array}{c}
\cos \theta(t) + w_x\\
\sin \theta(t) + w_y\\
u(t)/\rho
\end{array}
\right],\ \ \ \ \  u\in[-1,1]
\label{EQ:problem_wind}
\end{align}
where $\boldsymbol{w}:=(w_x,w_y)$ is the constant drift field induced by local winds/currents, and all other variables have been defined in Eq.~(\ref{Eq:problem1}). We consider a path planning problem of steering Eq.~(\ref{EQ:problem_wind}) from $\boldsymbol{z}_0$ to a position $(x_f,y_f)\in \mathbb{R}^2$ so that the resulting path is the shortest. To the authors's best knowledge, there is not a robust algorithm in the literature to solve this path planning problem efficiently. In fact, such a problem is equivalent to the MTIP if $\boldsymbol{v} = -\boldsymbol{w}$ \cite{Bakolas:2013}. Thus, we can use the algorithm in Section \ref{SE:Algorithm} to find the shortest Dubins path for any terminal point $(x_f,y_f)$ and any constant drift field $\boldsymbol{w}$. Four examples of the shortest Dubins paths in constant drift field are presented in Fig.~\ref{Fig:CaseE}.
 %where the constant drift field and the final point are indicated.

It is apparent to see from Fig.~\ref{Fig:CaseE} that none-straight curves along the shortest Dubins paths in constant drift fields are not circular arcs any more, which is consistent to the results in \cite{Techy:2009}.

\section{Conclusions}

While the MTIP is a fundamental problem in pursuit-evasion engagements, it is not exaggerate to say that, without any assumption on the distance between the initial position of the pursuer and the target's trajectory, this problem has not been well addressed. In this paper, through introducing three functions $F[E(t)]$, $L^+[E(t)]$, and $L^-[E(t)]$ and analyzing their continuity properties, it was shown that the solution of the MTIP lies in a sufficient family of 4 candidates. Moreover, the geometric properties of each candidate path was established, indicating that each candidate is a circular arc followed by either a circular arc or a straight line segment. When the target's velocity is constant, the geometric properties enabled formulating some nonlinear equations so that the length of each candidate was determined by a specific zero of the nonlinear equations. An efficient and robust algorithm was developed to find all the zeros of sufficiently smooth functions. As a result, the solution of the MTIP can be computed within a constant time. Since the MTIP with a constant target's velocity is equivalent to the RDP in a constant drift field, the developments of this paper also allowed  efficiently planning paths for aerial/marine vehicles in local winds/currents.

 \appendix

\section{Proofs for the theorems in Section \ref{SE:Syntheses}}
\label{Appendix:A}

%In this appendix, the proofs for Theorem \ref{TH:continuity} and Lemma \ref{LE:continuity} will be given.

Proof of Lemma \ref{TH:continuity}.
Let us consider that the target's trajectory $E(t)$ intersects the half circle $\mathcal{C}_r\cap \{y> 0\}$ at a time $\bar{t}> 0$. Without loss of generality, assume that the target enters into the circle $\mathcal{C}_r$ at $\bar{t}$ from outside. In such a case, there exists a small $\varepsilon > 0$ so that $E(t+ \varepsilon) \in \mathrm{int} (\mathcal{D}_r)$ and $E(t- \varepsilon) \not \in  \mathcal{D}_r$ where the notation $\mathrm{int}(\cdot)$ denotes the interior of a set.
\begin{figure}[h]
	\centering
	\includegraphics[width = 6cm]{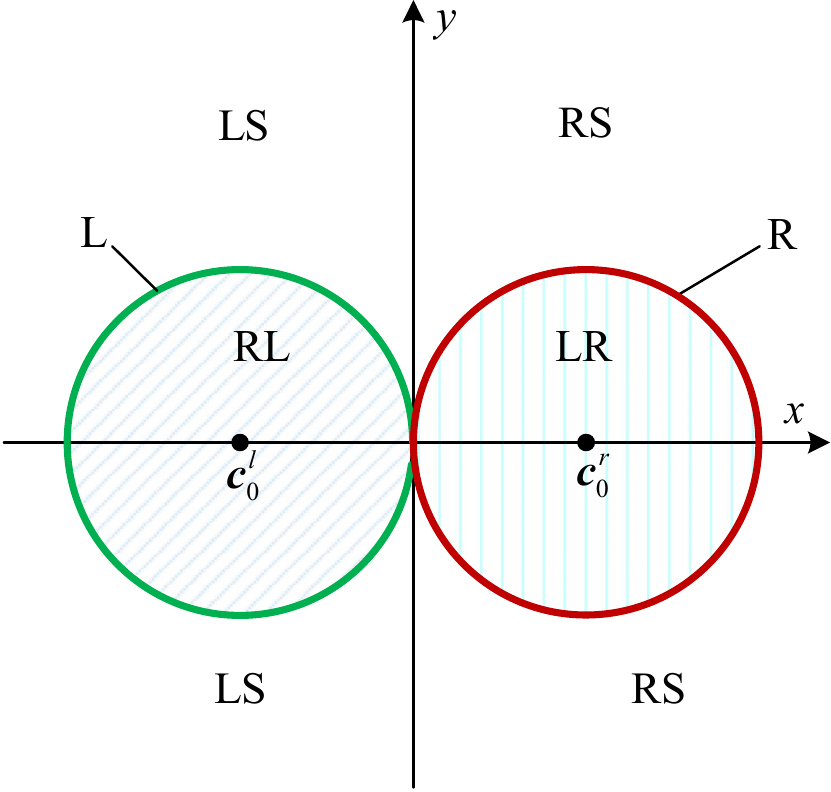}
	\caption{The regions for which the path of RDP terminates by different types. \cite{Boissonnat:1994}}
	\label{Fig:Types_Region}
\end{figure}
According to \cite{Boissonnat:1994},  the shortest paths of RDP from $\boldsymbol{z}_0$ to $E(\bar{t}-\varepsilon)$, $E(\bar{t})$,  and $E(\bar{t}+\varepsilon)$ are of types RS$_{d(\varepsilon)}$, R, and LR$_{\alpha(\varepsilon)}$, respectively, as shown by Fig.~\ref{Fig:Types_Region}.
\begin{figure}[h]
	\centering
	\includegraphics[width = 7cm]{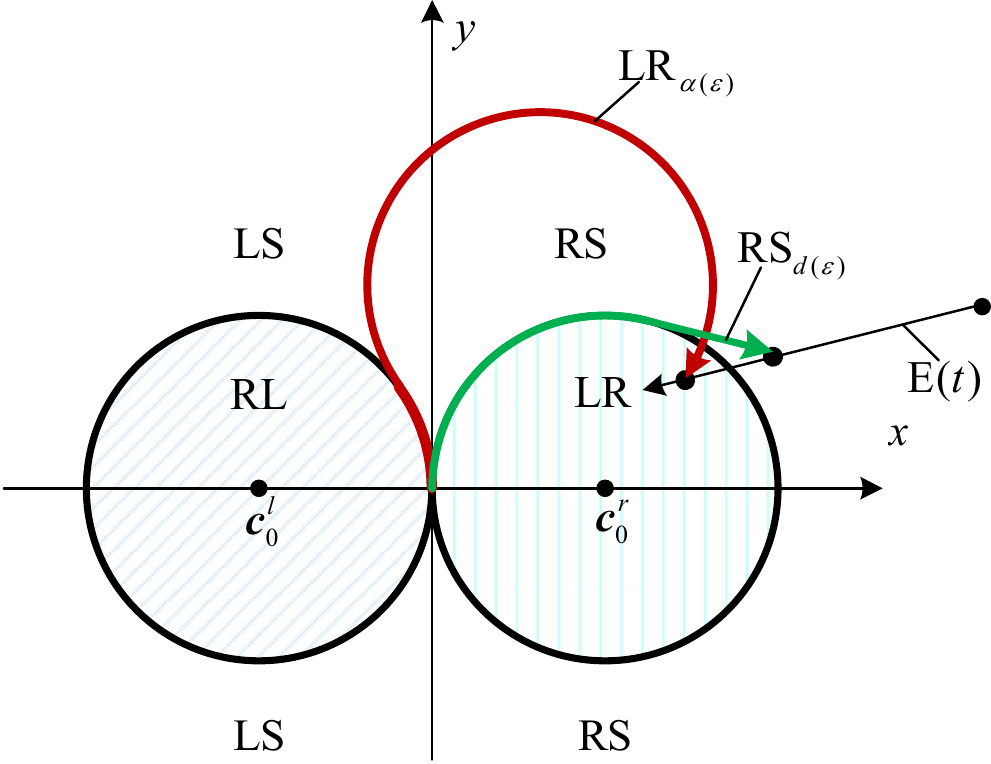}
	\caption{The geometry for the path of RDP with terminal points $E(\bar{t}+\varepsilon)$, $E(\bar{t})$, and $E(\bar{t}-\varepsilon)$.}
	\label{Fig:Intersection}
\end{figure}
It is also known from \cite{Boissonnat:1994} that  for any sufficiently small $\varepsilon>0$ we have $\alpha(\varepsilon)>\pi$, indicating that the function $F[E(t)]$ is discontinuous when the target moves from the circle $\mathcal{C}_r$ into its inside. For the case that the target's trajectory $E(t)$ intersects the left half circle $\mathcal{C}_l \cap \{y > 0\}$, it can be proven in the same way that the function $F[E(t)]$ is discontinuous at the intersection point.

% However,  As $\varepsilon$ increases, the type of the solution path of RDP changes from $L_{\alpha(\varepsilon)}R$ to $R$.
%It should be noted that, with the increase of $t$ from $\bar{t}-\varepsilon$ to $\bar{t}+\varepsilon$, the  path of RDP changes its type from RS, via R, to LR. According to \cite[Section 3.3]{Boissonnat:1994}, we have
%$\underset{\varepsilon \rightarrow 0}{\lim}\  d(\varepsilon) = 0$
%but
%$\underset{\varepsilon \rightarrow 0}{\lim}\   \alpha(\varepsilon) \neq 0$. {\color{blue}****can use that $\alpha>pi+pi/3$**** }Hence, the function $F[E(t)]$ is not continuous at $\bar{t}$.

From now on, we prove the necessity that the function $F[E(\bar{t})]$ is continuous if $E(\bar{t})\not \in \{(x,y)\in \mathcal{C}_r\cup \mathcal{C}_l \big | y> 0\}$.  Let us choose a time $\bar{t}>0$ so that $E(\bar{t})\in \{(x,y)\in \mathbb{R}^2 \big | (x,y)\not \in \mathcal{D}_r\ \text{and}\ x>0\}$.  According to \cite[Section 3.3]{Boissonnat:1994}, given a sufficiently small $\varepsilon > 0$, the path of RDP from $\boldsymbol{z}_0$ to $E(t+\delta )$ for any $\delta \in [-\varepsilon ,\varepsilon]$ is of type RS. Thus, the function $F[E(t)]$ is continuous at $\bar{t}$ if $E(\bar{t}) \in \{(x,y)\in \mathbb{R}^2 \big | (x,y)\not \in \mathcal{D}_r\ \text{and}\ x>0\}$.

We choose a time $\bar{t}>0$ so that $E(\bar{t})\in \mathrm{int}(\mathcal{D}_r)$.  According to \cite[Section 3.3]{Boissonnat:1994}, given a sufficiently small $\varepsilon > 0$, the path of RDP from $\boldsymbol{z}_0$ to $E(t+\delta )$ for any $\delta \in [-\varepsilon ,\varepsilon]$ keeps being of type LR.  Thus, the function $F[E(t)]$ is continuous at $\bar{t}$ if $E(\bar{t}) \in \mathrm{int}(\mathcal{D}_r)$.

Then, we consider the rest case that there exists a time $\bar{t}> 0 $ so that $E(\bar{t}) \in \mathcal{C}_r \cap \{y \leq 0\}$. Without loss of generality, let us assume that the target enters into the circle $\mathcal{C}_r$ at $\bar{t}$. Then, there exists a small $\varepsilon > 0$ so that $E(t+\varepsilon) \in \mathrm{int}(\mathcal{D}_r)$ and $E(t- \varepsilon )\not \in \mathrm{int}(\mathcal{D}_r)$.  According to \cite[Section 3.3]{Boissonnat:1994}, the paths of RDP from $\boldsymbol{z}_0$ to $E(\bar{t}-\varepsilon )$, $E(\bar{t})$, and $E(\bar{t}+\varepsilon)$ are of types RS, R, and LR. Note that if $\varepsilon$ approaches to zero, the left turning circle $L$ and the straight line segment vanish. Thus, the function $F[E(t)]$ is continuous at $\bar{t}$ if $E(\bar{t}) \in \mathcal{C}_r \cap \{y \leq 0\}$.

If the time $\bar{t}>0$ is chosen so that $E(\bar{t}) \in \{x< 0\}$, the continuity of $F[E(t)]$ can be proven in the same way, completing the proof.
 $\square$

\begin{lemma}\label{RE:1}
The minimum time $t_m > 0$ for the pursuer to intercept the target is not smaller than $F[E(t_m)]$, i.e., $t_m \geq F[E(t_m)]$.
\end{lemma}
Proof of Lemma \ref{RE:1}. According to the definitions of $F[\cdot]$ and $E(\cdot)$, for any $t\geq 0$ the value $F[E(t)]$ denotes the minimum time for the Dubins vehicle to move from $\boldsymbol{z}_0$ to $E(t)$. Thus, the duration for the Dubins vehicle to move from $\boldsymbol{z}_0$ to $E(t_m)$ must be greater than $F[E(t_m)]$, completing the proof.
$\square$

\iffalse

Proof of Lemma \ref{LE:continuous}. Consider a new function $\hat{F}(t) =t - F[E(t)]$, and we first prove that $t_m$ is a zero of $\hat{F}(t)$, i.e., $\hat{F}(t_m) = 0$.  According to Lemma \ref{RE:1}, we must have $\hat{F}(t_m ) \geq 0 $. In order to prove that $\hat{F}(t_m) = 0$ holds, we  by contradiction assume that $\hat{F}(t_m ) > 0 $. Note that $\hat{F}(0) < 0$ and $\hat{F}(t)$ is continuous. Thus, according to intermediate value theorem, there exists a time $\bar{t}\in (0,t_m)$ so that $\hat{F}(\bar{t}) = 0$, indicating $\bar{t} = F[E(\bar{t})]$. Thus, the target can be intercepted by the pursuer at a time $\bar{t}$ smaller than $t_m$. By contraposition, we have $t_m = F[E(t_m)]$.  In case that the function $F[E(\cdot)]$ has multiple fixed points,  the minimum interception time $t_m>0$ must be the minimum fixed point of $F[E(\cdot)]$, concluding the proof.  $\square$

\fi

Proof of Lemma \ref{LE:fixed_RDP}. By contradiction, assume that the solution path of the MTIP does not follow that of the RDP from $\boldsymbol{z}_0$ to $E(t_m)$. Notice that $F[E(t_m)]$ denotes the time for the Dubins vehicle to follow the solution path of RDP from $\boldsymbol{z}_0$ to $E(t_m)$. Since the solution of RDP is the optimal path from $\boldsymbol{z}_0$ to $E(t_m)$, the contradicting assumption implies $t_m \neq F[E(t_m)]$. This contradicts the assumption of the lemma, completing the proof.
$\square$

Proof of Theorem \ref{LE:Occurance_R1}. (1) By contradiction, assume that $t_m$ is not a fixed point of $F[E(\cdot)]$. According to Lemma \ref{RE:1}, this contradicting assumption implies  $t_m > F[E(t_m)]$. In view of Lemma \ref{TH:continuity}, if $E(t_m) \in \mathrm{int}(\mathcal{R}_1\cup \mathcal{R}_2)$, we have that $F[E(\cdot)]$ is continuous around $t_m$. Thus, there exists a sufficiently small $\varepsilon>0$ so that
\begin{align}
t_m - \varepsilon > F[E(t_m - \varepsilon)].
\label{EQ:theorem_R1}
\end{align}
Since $E(t_m)\in \mathrm{int}(\mathcal{R}_1\cup \mathcal{R}_2)$ and $E(t)$ is continuous, if $\varepsilon > 0$ is sufficiently small, it holds that $E(t_m - \varepsilon)\in  \mathrm{int}(\mathcal{R}_1\cup \mathcal{R}_2)$. Then, according to \cite[Theorem 2]{Ding:2019}, for any $t> F[E(t_m - \varepsilon)]$ there exists a feasible Dubins path with a duration of $t$ from initial condition $\boldsymbol{z}_0$ to the point $E(t_m-\varepsilon)$. For this reason, since $t_m - \varepsilon > F(t_m - \varepsilon )$ by Eq.~\eqref{EQ:theorem_R1}, it follows that there exists a feasible Dubins path with a duration of $t_m - \varepsilon$ to reach the point $E(t_m - \varepsilon)$. This means that the interception between the pursuer and the target can occur at a time $t_m - \varepsilon$, which is smaller than $t_m$. By contraposition, the proof of the first statement is completed.

(2) Combining the first statement and Lemma \ref{LE:fixed_RDP}, the second statement holds apparently.
$\square$

Proof of Lemma \ref{LE:continuity}. Given a path of type L$_u$R$_v$, the terminal point is expressed by \cite{Boissonnat:1994}
\begin{align}
(x,y) = \left(
\begin{array}{l}
-\rho + 2\rho \cos u +  \rho \cos (u+ \pi - v)\\
2\rho \sin u + \rho \sin (u+ \pi - v)
\end{array}\right)
\nonumber
\end{align}
According to this formula, if the terminal point continuously changes, the values of $u$ and $v$ continuously change as well. Thus, the functions $L^-[E(t)]$ and $^+[E(t)]$ are continuous at $\bar{t}$, completing the proof.
$\square$

Proof of Theorem \ref{LE:R3}.  In view of the second statement of Lemma \ref{LE:Ding}, for any given $t$ in the open interval$ (L^-[E(t_m)],L^+[E(t_m)])$,  any feasible Dubins path with a duration of $t$ cannot start from $\boldsymbol{z}_0$ and terminate at $E(t_m)$. Hence,  in order to prove this theorem, we just need to prove that $t_m$ lies on the boundary of the interval $[F[E(t_m)],L^-[E(t_m)]]\cup [L^+[E(t_m)],+\infty)$.

By contradiction, assume $t_m \in (F[E(t_m)],L^-[E(t_m)])$. Note that $F[E(t)]$  and $L^-[E(t)]$ are continuous at $t_m$ according to Lemma \ref{TH:continuity} and Lemma \ref{LE:continuity}, respectively. Hence, there exists a sufficiently small $\varepsilon > 0$ so that $F[E(t_m - \varepsilon )] < t_m - \varepsilon < L^-[E(t_m - \varepsilon)]$. According to the third statement of Lemma \ref{LE:Ding},  for any $t$ in the interval $[F[E(t_m - \varepsilon)],L^-[E(t_m - \varepsilon)]]$, there exists a feasible Dubins path with a duration of $t$ from $\boldsymbol{z}_0$ to the point $E(t_m - \varepsilon)$. Thus, the target can be intercepted by the pursuer  at the point $E(t_m - \varepsilon)$ by moving along a feasible Dubins path with a duration of $t_m - \varepsilon$, smaller than $t_m$. By contraposition, we have that $t_m$ does not lie in the interior of the interval $[F[E(t_m)],L^-[E(t_m)]]$, i.e., $t_m \not\in (F[E(t_m)],L^-[E(t_m)])$.

From now on, we assume that $t_m$ lies in the open interval $ (L^+[E(t_m)],+\infty)$ by contradiction.  Note that $L^+[E(t)]$ is continuous at $t_m$ according to Lemma \ref{LE:continuity}. Therefore, there exists a sufficiently small $\varepsilon > 0$ so that $L^+(E(t_m - \varepsilon)) < t_m - \varepsilon$. Since $E(t_m)\in \mathrm{int}(\mathcal{R}_3\cap\{x>0\})$, if $\varepsilon>0$ is sufficiently small, it holds $E(t_m-\varepsilon)\in \mathrm{int}(\mathcal{R}_3\cap\{x>0\})$.  Then, based on the third statement of Lemma \ref{LE:Ding}, for any $t$ in the semi-open interval $[L^+(E(t_m-\varepsilon)),+\infty)$, there exists a feasible Dubins path from $\boldsymbol{z}_0$ to the point $E(t_m - \varepsilon)$ with a duration of $t$, smaller than $t_m$. By contraposition, we have that $t_m$ does not lie in the open interval $(L^+[E(t_m)],+\infty)$.

To sum up, it is concluded that $t_m$ lies on the boundary of $[F[E(t_m)],L^-[E(t_m)]]\cup [L^+[E(t_m)],+\infty)$, indicating that $t_m \in \{t> 0 | t = F[E(t)], L^-[E(t)], or  L^+[E(t)]\}$.  Since $t_m$ is the minimum intercept time, we have that Eq.~\eqref{EQ:tm_R3} holds, completing the proof.
$\square$

\section{Proofs for the lemmas in Section \ref{SE:Algorithm}}\label{Appendix:B}

Proof of Lemma \ref{LE:extreme_to_zero}. This lemma is a direct result of the intermediate value theorem. $\square$

\begin{figure}[h]
	\centering
	\includegraphics[width =8cm]{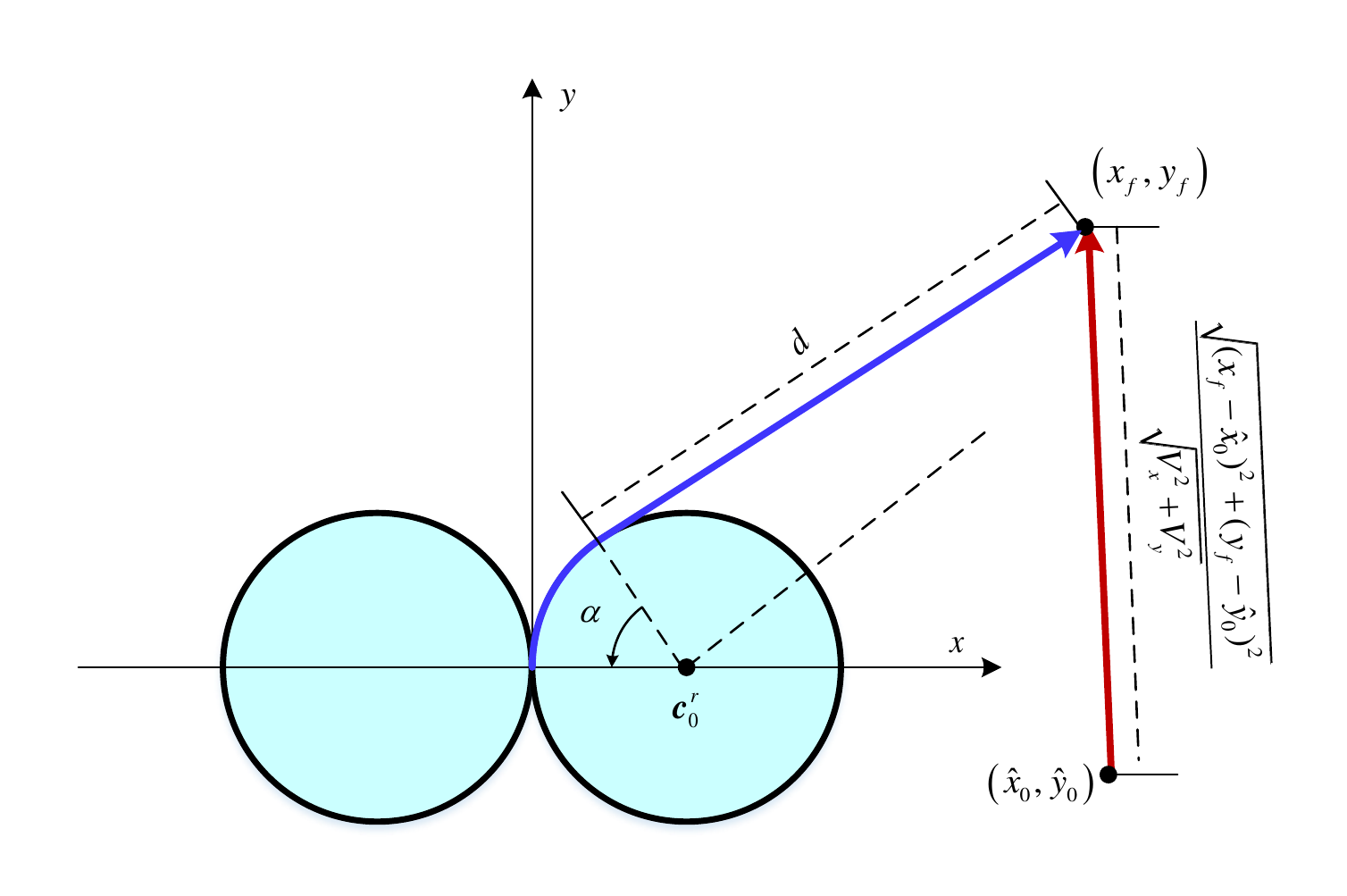}
	\caption{The geometry of the solution path of RDP with a type of RS.}
	\label{Fig:LS1}
\end{figure}

Proof of Lemma \ref{LE:analytic_CS}.
Let $\alpha \in [0,2\pi]$ be the radian of the right-turning circular arc R in the path of type RS, as presented in Fig.~\ref{Fig:LS1}. By geometric analysis, we have
%\begin{small}
\begin{align}
%\begin{cases}
\boldsymbol{c}_0^r + \rho \left[
\begin{array}{c}
 -\cos (\alpha )\\
 \sin (\alpha )
\end{array}
\right] + d\left[
\begin{array}{c}
\cos(\frac{\pi}{2}- \alpha)\\
\sin (\frac{\pi}{2} - \alpha)
\end{array}
\right]
=\left[
\begin{array}{c}
x_f\\
y_f
\end{array}
\right]
%\end{cases}
\label{LS0}
\end{align}
%\end{small}
where $d>0$ is the length of the straight line segment and $(x_f,y_f)$ is the intercept point.
 Because the pursuer and the target arrive at $(x_f,y_f)$ simultaneously, we have
\begin{align}
\rho \alpha + d = \frac{\sqrt{(x_f - \hat{x}_0)^2 + (y_f - \hat{y}_0)^2}}{\sqrt{v_x^2 + v_y^2}}\label{LS1}
\end{align}
%where $n=0$ if $\alpha \in [0,\pi/2]$ and $n=1$ if $\alpha\in (\pi/2,2\pi)$.
Note that we have
\begin{align}
\frac{y_f-\hat{y}_0}{x_f-\hat{x}_0}=\frac{v_y}{v_x}\nonumber
\end{align}
indicating
\begin{align}
y_f=\frac{v_y}{v_x}x_f-\frac{v_y}{v_x} \hat{x}_0+\hat{y}_0.
\label{LS2}
\end{align}
Combining Eq.~(\ref{LS1}) and Eq.~(\ref{LS2})  yields
\begin{align}
d=\frac{x_f-\hat{x}_0}{v_x}-\rho \alpha
\label{LS3}
\end{align}
Substituting Eq.~(\ref{LS2}) and Eq.~(\ref{LS3}) into Eq. ~\eqref{LS0} to eliminate $x_f$ and $d$, we have
%{\color{blue}
\begin{align}
\frac{a_1\!+\! a_2\sin \alpha\!-\!\rho \cos \alpha\!-\!\rho \alpha \sin \alpha}{a_3\!+\!a_2\cos \alpha\!+\!\rho \sin \alpha\!-\!\rho \alpha \cos \alpha} =\frac{v_x-\sin \alpha}{v_y-\cos \alpha}
\label{LS4}
\end{align}
%}
where
%{\color{blue}
\begin{align}
\begin{cases}
a_1=\rho\\
a_2=\frac{\hat{x}_0}{v_x}\\
a_3=-\hat{y}_0 + \frac{v_y}{v_x}\hat{x}_0
\end{cases}\nonumber
\end{align}
%}
Rearranging Eq.~(\ref{LS4}), we  get
\begin{align}
A_1 \sin \alpha +A_2 \cos \alpha + \alpha(A_3 \cos \alpha +A_4  \sin \alpha) +A_5=0
\label{LS5}
\end{align}
where
%{\color{blue}
\begin{align}
\begin{cases}
A_1=a_2 v_y + a_3 - \rho v_x\\
A_2=- a_1 - \rho v_y - a_2 v_x\\
A_3=\rho v_x \\
A_4=- \rho v_y \\
A_5=a_1 v_y - a_3 v_x + \rho
\end{cases}\nonumber
\end{align}
%}
This concludes the proof of Lemma  \ref{LE:analytic_CS}. $\square$

\begin{figure}[h!]
	\centering
	\includegraphics[width = 8cm]{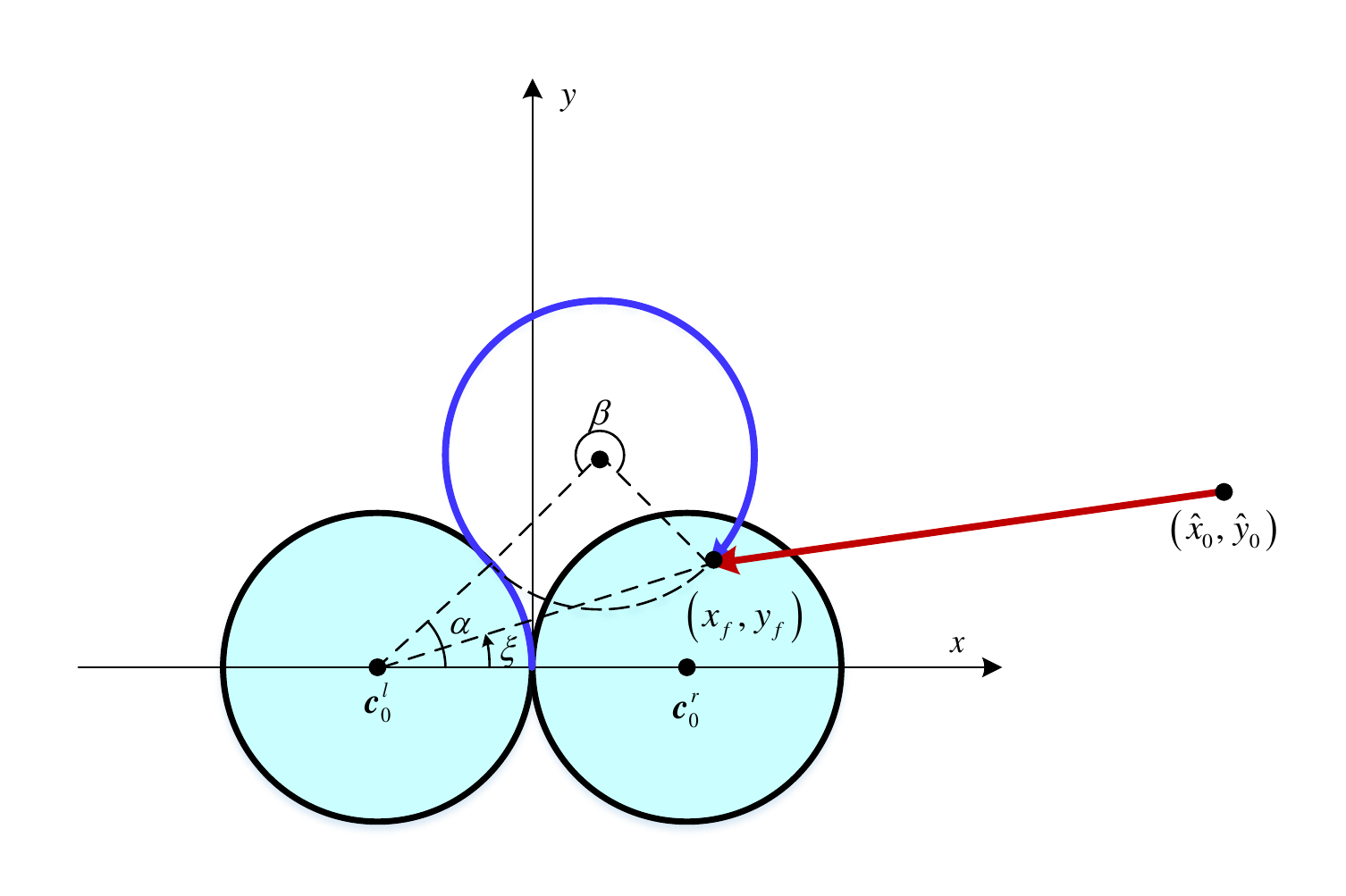}
	\caption{The geometry of the solution path of RDP with a type of LR.}
	\label{Fig:RL1}
\end{figure}

Proof of Lemma \ref{LE:analytic_CC}. Let $\alpha\in[0,2\pi]$ and $\beta\in[\pi,2\pi]$ be the radians of $L$ and $R$, respectively, as presented in Fig.~\ref{Fig:RL1}.

Note that $$\boldsymbol{c}_0^l+2\rho\left[
\begin{array}{c}
\cos(\alpha)\\
\sin( \alpha)
\end{array}
\right] $$
and
$$
\left[
\begin{array}{c}
x_f\\
y_f
\end{array}
\right]
+ \rho
\left[
\begin{array}{c}
\cos(\alpha - \beta)\\
\sin( \alpha - \beta )
\end{array}
\right]$$
are located at the same point (the center of the circle coinciding with R). Thus, we have
\begin{align}
 \boldsymbol{c}_0^l + 2\rho
\left[
\begin{array}{c}
\cos(\alpha)\\
\sin(\alpha )
\end{array}
\right]
=
\left[
\begin{array}{c}
x_f\\
y_f
\end{array}
\right]
+ \rho
\left[
\begin{array}{c}
\cos( \alpha - \beta )\\
\sin( \alpha - \beta)
\end{array}
\right]
\label{RL2}
\end{align}

As the pursuer and the target arrive at $(x_f,y_f)$  simultaneously, it follows that
\begin{align}
\rho(\alpha + \beta) =
\frac{\sqrt{(x_f - \hat{x}_0)^2 + (y_f - \hat{y}_0)^2}}{\sqrt{v_x^2 + v_y^2}}
\label{RL0}
\end{align}
We also have
\begin{align}
\frac{y_f-\hat{y}_0}{x_f-\hat{x}_0}=\frac{v_y}{v_x}
\label{RL1}
\end{align}
Combining Eq.~\eqref{RL0} with Eq.~\eqref{RL1} leads to
\begin{align}
\begin{split}
 \alpha + \beta=\frac{x_f-\hat{x}_0}{\rho v_x }= \frac{y_f-\hat{x}_0}{\rho vp_y}
 \end{split}
\label{RL11}
\end{align}
By the law of cosines, we also have
\begin{align}
4 \rho^2+\rho^2-4\rho^2 \cos (2\pi-\beta)=s^2
\label{RL4}
\end{align}
which leads to
\begin{align}
\begin{cases}
\cos \beta =
\frac{5\rho^2-s^2}{4 \rho^2}\\
\sin \beta = \pm \sqrt{1 - \left[\frac{5\rho^2-s^2}{4 \rho^2}\right]^2}
\end{cases}
\label{RL5}
\end{align}
Denote by $\xi>0$ the angle between $x$-axis and the vector from $\boldsymbol{c}_0^l$ to $(x_f,y_f)$. Then, by using the law of sines, we have
\begin{align}
\frac{\sin (2\pi - \beta)}{s}=\frac{\sin (\alpha - \xi)}{\rho}
\label{RL3}
\end{align}
where $s = \| \boldsymbol{c}_0^l - (x_f,y_f)\|$  is the Euclidean distance between $\boldsymbol{c}_0^l$ and the interception point $(x_f,y_f)$.
Taking into account the expression of $s$, from Eq.~(\ref{RL2}) we get
\begin{align}
0  = s^2  -  4 \rho ( x_f + \rho) \cos \alpha  - 4 \rho  y_f \sin \alpha + 3 \rho^2
\label{RL7}
\end{align}
Note that $\cos \xi=\frac{x_f + \rho}{s}$ and $\sin  \xi=\frac{y_f}{s}$. Then, according to Eq.~(\ref{RL3}), we have
\begin{align}
\sin \beta = \frac{y_f}{\rho} \cos \alpha - \frac{x_f  + \rho}{\rho} \sin \alpha
\label{RL9}
\end{align}
Substituting Eq.~(\ref{RL5}) and Eq.~(\ref{RL7}) into Eq.~(\ref{RL9}), we can get
\begin{align}
\begin{cases}
\sin \alpha =  \pm \frac{\rho ( x_f + \rho)\sqrt{1\!-\![\frac{s^2-5\rho^2}{4\rho^2}]^2}}{s^2} +y_f \frac{3\rho^2+s^2}{4\rho s^2}\\
\cos \alpha\!= \pm \!\frac{\rho y_f \sqrt{1\!-\![\frac{s^2-5\rho^2}{4\rho^2}]^2}}{s^2}\!+(x_ f + \rho)\frac{3\rho^2+s^2}{4\rho s^2}
\end{cases}
\label{RL10}
\end{align}

Set $\eta = \alpha + \beta$. Then, we have
\begin{align}
\begin{split}
\sin \eta = \sin \alpha \cos\beta + \cos \alpha \sin \beta\\
\cos \eta = \cos \alpha \cos \beta - \sin \alpha \sin \beta
\end{split}
\label{EQ:sin_cos_eta}
\end{align}
Substituting Eq.~(\ref{RL5}) and Eq.~(\ref{RL10}) into Eq.~\eqref{EQ:sin_cos_eta}, we have
\begin{equation}
\left\{
\begin{array}{l}
\pm ( x_f +\rho )\frac{s^2-\rho^2}{2\rho s^2}\sqrt{1-[\frac{s^2-5\rho^2}{4\rho^2}]^2}=\sin \eta
+\frac{(y_0^r - y_f )(s^4-6\rho^2s^2-3\rho^4)}{8\rho^3s^2}\\
\specialrule{0em}{0.5ex}{0.5ex}
\pm y_f\frac{s^2-\rho^2}{2\rho s^2}\sqrt{1-[\frac{s^2-5\rho^2}{4\rho^2}]^2}
 =\cos\eta+\frac{(x_0^r - x_f)(s^4-6\rho^2s^2-3\rho^4)}{8\rho^3s^2} \\
\end{array}\right.
\label{RL13}
\end{equation}
where $(x_0^r,y_0^r) = \boldsymbol{c}_0^r$.
Rewriting Eq.~(\ref{RL13}) yields
\begin{align}
y_f \sin\eta& + ( x_f +\rho) \cos \eta+\frac{s^4-6 \rho^2 s^2 - 3\rho^4}{8\rho^3} =0
\label{RL14}
\end{align}
According to Eq.~(\ref{RL11}), we have
\begin{align}
\begin{cases}
x_f=\eta v_x \rho +\hat{x}_0\\
y_f=\eta v_y \rho +\hat{y}_0
\end{cases}
\label{RL15}
\end{align}
Substituting  Eq.~(\ref{RL15}) into  Eq.~(\ref{RL14}), we eventually have
\begin{align}
F(\eta)& \overset{\triangle}{=}B_1\eta^4 + B_2\eta^3 +  B_3\eta^2 + B_4\eta + B_5 +B_6\cos \eta +  B_7 \sin \eta  + \eta(B_8  \cos \eta  + B_9  \sin \eta)= 0
\label{RL18}
\end{align}
where
\begin{align}
\begin{split}
B_1&=C_a^2\\
B_2&=2C_aC_b\\
B_3&=C_b^2+2C_aC_c-6\rho^2C_a\\
B_4&=2C_cC_b-6\rho^2C_b\\
B_5&=C_c^2-6\rho^2C_c-3\rho^4\\
B_6&=8\rho^3(\rho + \hat{x}_0)\\
B_7&=8\rho^3\hat{y}_0\\
B_8&=8\rho^4V_x\\
B_9&=8\rho^4V_y
\end{split}\nonumber
%\label{RL19}
\end{align}
with
\begin{align}
\begin{split}
C_a &=V_x^2\rho^2+V_y^2\rho^2\\
C_b& =+2(\rho + \hat{x}_0)V_x\rho+2\hat{y}_0V_y\rho\\
C_c&=(\rho + \hat{x}_0)^2+\hat{y}_0^2
\end{split}\nonumber
%\label{RL20}
\end{align}

\bibliographystyle{unsrt}  
\bibliography{references}  %%% Remove comment to use the external .bib file (using bibtex).

%\begin{thebibliography}{1}

%\bibitem{kour2014real}
%George Kour and Raid Saabne.
%\newblock Real-time segmentation of on-line handwritten arabic script.
%\newblock In {\em Frontiers in Handwriting Recognition (ICFHR), 2014 14th
%  International Conference on}, pages 417--422. IEEE, 2014.

%\bibitem{kour2014fast}
%George Kour and Raid Saabne.
%\newblock Fast classification of handwritten on-line arabic characters.
%\newblock In {\em Soft Computing and Pattern Recognition (SoCPaR), 2014 6th
 % International Conference of}, pages 312--318. IEEE, 2014.

%\bibitem{hadash2018estimate}
%Guy Hadash, Einat Kermany, Boaz Carmeli, Ofer Lavi, George Kour, and Alon
%  Jacovi.
%\newblock Estimate and replace: A novel approach to integrating deep neural
%  networks with existing applications.
%\newblock {\em arXiv preprint arXiv:1804.09028}, 2018.

%\end{thebibliography}

\end{document}